%% file: manuscript.tex
\documentclass{article}

\PassOptionsToPackage{numbers,sort&compress}{natbib}
\usepackage[preprint]{neurips_2026}
\newcommand{\keywords}[1]{\par\noindent\textbf{Keywords}\quad #1\par}
\usepackage[utf8]{inputenc}
\usepackage[T1]{fontenc}
\usepackage{array}
\usepackage{microtype}
\usepackage{url}
\usepackage{doi}

\usepackage{amsmath,amssymb,amsthm}
\usepackage{dsfont}
\usepackage{bm}
\usepackage{graphicx}
\usepackage{booktabs}
\usepackage{threeparttable}
\usepackage{placeins}
\usepackage{flafter}
\usepackage{needspace}
\usepackage{siunitx}
\usepackage{algorithm}
\usepackage{algpseudocode}
\usepackage{tikz}
\usetikzlibrary{arrows.meta,positioning,shapes.geometric}

\newcommand{\indicator}[1]{\mathds{1}_{#1}}

\usepackage{hyperref}
\hypersetup{hypertexnames=false,hidelinks,pdfdisplaydoctitle=true}

\title{Reachability-Informed Reinforcement Learning for Multi-Impulse Interplanetary Transfers}
\author{Yashdeep Chaudhary
\thanks{%
\begin{tabular}[t]{@{}l@{}}
Corresponding author.\\
Email address:
\href{mailto:yashdeep.chaudhary@auckland.ac.nz}
{\texttt{yashdeep.chaudhary@auckland.ac.nz}}
(Yashdeep Chaudhary)
\end{tabular}%
}\\
\normalfont\small Faculty of Engineering\\
\normalfont\small Waipapa Taumata Rau--\\
\normalfont\small University of Auckland\\
\normalfont\small Auckland 1010, New Zealand
\And
Roberto Armellin\\
\normalfont\small Faculty of Engineering\\
\normalfont\small Waipapa Taumata Rau--\\
\normalfont\small University of Auckland\\
\normalfont\small Auckland 1010, New Zealand
\And
Harry Holt\\
\normalfont\small Advanced Concepts Team\\
\normalfont\small ESTEC, European Space Agency\\
\normalfont\small Noordwijk 2201 AZ\\
\normalfont\small The Netherlands}
\hypersetup{pdftitle={Reachability-Informed Reinforcement Learning for Multi-Impulse Interplanetary Transfers},pdfauthor={Yashdeep Chaudhary, Roberto Armellin, Harry Holt}}

\makeatletter
\renewcommand{\@notice}{}
\makeatother
\usepackage[font=small,labelfont=bf]{caption}
\usepackage{bookmark}
\hypersetup{bookmarksnumbered=true,bookmarksopen=false}
\makeatletter
\renewcommand{\@toptitlebar}{\hrule height 2pt\vskip 0.25in\vskip -\parskip}
\renewcommand{\@bottomtitlebar}{\vskip 0.29in\vskip -\parskip\hrule height 2pt\vskip 0.09in}
\makeatother

\AtBeginDocument{\newgeometry{textwidth=6.5in,textheight=9in,top=1in,headheight=12pt,headsep=25pt,footskip=30pt}}

\begin{document}
\maketitle

\begin{abstract}
Reinforcement learning offers the prospect of a reusable sequential decision-making mechanism for spacecraft trajectory design, motivating policy interfaces that connect learned decisions to the underlying maneuver geometry. This paper develops Reachability Analysis-Informed Reinforcement Learning (RARL) for deterministic multi-impulse interplanetary transfers, placing intermediate waypoint selection at the center of the learned decision process. Local first-order reachability maps bounded velocity perturbations into an ellipsoidal set of next-node positions, within which the policy selects its waypoint. Lambert reconstruction then determines the corresponding maneuver to reach this selected waypoint along a dynamically consistent ballistic arc, coupling learned transfer-geometry selection with classical astrodynamics. A terminal two-impulse reconstruction completes the rendezvous, supported by a linear maneuver-demand assessment used for reward shaping. Numerical studies characterize this interface on a two-body Earth--Mars benchmark. Across three independent training runs, RARL achieves a mean maneuver cost of 10.23 km/s, 1.72\% above a validated local sequential convex programming reference. Training over dispersed initial states extends policy reuse across a departure family with fixed target state and transfer duration. Each of the three independently trained multi-state policies completes all 10,000 held-out Monte Carlo departures without impulse-cap violations, compared with a mean feasibility rate of 6.49\% for single-state policies. This broader sampled feasibility is accompanied by a 0.61\% increase in mean nominal maneuver cost, without further training across departures. These results demonstrate that a reachability-informed decision interface supports benchmark-quality trajectory construction and policy reuse across dispersed departure conditions.
\end{abstract}

\keywords{
Interplanetary transfer, Lambert problem, reachability analysis,
reinforcement learning, trajectory\linebreak[4]
\mbox{optimization}, spacecraft guidance, multi-impulse transfers.
}

\input{sections/introduction.tex}
\input{sections/problem-formulation.tex}
\input{sections/methodology.tex}
\input{sections/validation.tex}
\input{sections/case-study-results.tex}
\FloatBarrier
\input{sections/conclusions.tex}
\input{sections/app-secondary-proxy.tex}

\FloatBarrier

\end{document}

%% file: sections/introduction.tex
\section{Introduction}
\label{sec:introduction}

Even in deterministic settings, interplanetary transfers are shaped by nonlinear dynamics, limited maneuver authority, and stringent terminal requirements. In a multi-impulse rendezvous, these requirements couple the allocation of individual maneuvers to the geometry and cost of the complete transfer. Classical indirect optimal-control methods and direct transcription provide established tools for computing such trajectories~\cite{Betts1998-ey,Conway2010-td}. Their development also extends to joint trajectory and feedback design under uncertainty~\cite{benedikter2022ConvexApproachCovariance}. Alongside the design of an individual reference trajectory, mission analysis often requires the assessment of a family of departure conditions. This motivates complementing trajectory optimization with reusable decision mappings that can generate transfers from different initial states without solving the complete optimization problem anew for each departure.

Reinforcement learning (RL) offers a means of constructing these mappings by learning a state-to-action policy through simulated interaction~\cite{sutton2018ReinforcementLearningIntroduction}. The main attraction is the opportunity to perform policy training offline and reuse the resulting mapping across the conditions for which it has been trained and assessed. Within the broader application of artificial intelligence to spacecraft guidance, dynamics, and control~\cite{Izzo2019AISurvey}, this approach has been explored for interplanetary trajectory design under uncertainty~\cite{zavoli2021ReinforcementLearningRobust,federici2024RobustInterplanetaryTrajectory} and autonomous cislunar orbit transfers~\cite{federici2021AutonomousGuidanceCislunar}. These applications establish a basis for investigating how learned policies can accommodate variations in transfer conditions while meeting maneuver and terminal requirements.

A complementary line of work embeds learning within an established guidance structure. Reinforced Q-laws~\cite{holt2021OptimalQlawsReinforcement} and Lyapunov controllers for lunar transfers~\cite{Holt2024RLC} use learned parameters to adapt the underlying guidance law. Reinforced Lyapunov guidance has also been combined with precomputed ballistic capture trajectories for low-thrust cislunar transfers~\cite{chaudhary2024LowThrustCisLunarTransfers}. More recently, chance-constrained RL has been used to augment deterministic reference trajectories with learned feedforward and feedback corrections under sampled uncertainty~\cite{chaudhary2026DistributionAgnosticRCCRL}. These studies illustrate how learning can be assigned a defined role within a physically structured trajectory-design framework. For multi-impulse transfers, a related opportunity is to incorporate the relationship between maneuver authority and attainable transfer geometry directly into the policy's decision space.

This relationship is particularly relevant to rendezvous, where each maneuver changes both the subsequent ballistic arc and the maneuver demand required to meet the terminal state. Classical sensitivity maps and boundary-value solutions already describe aspects of this relationship~\cite{battin2001IntroductionMathematicsMethods,vallado2013FundamentalsAstrodynamicsApplications}. They therefore provide a basis for expressing a policy decision in terms of an intermediate target and reconstructing the maneuver required to reach it. The resulting design question is how to give the policy a geometrically meaningful set of next-node targets that reflects the available impulsive authority.

Reachability analysis provides a natural basis for this construction by characterizing attainable states under specified dynamics and control bounds. Hamilton--Jacobi formulations describe reachable sets for nonlinear systems and dynamic games~\cite{mitchell2005time,bansal2017hamilton}, while ellipsoidal and zonotopic methods provide geometric representations and approximations of such sets~\cite{kurzhanski2000ellipsoidal,althoff2010computing}. In astrodynamics, reachable domains have been studied for single-impulse maneuvers~\cite{xue2012ReachableDomainSpacecraft}, tangent impulses subject to trajectory-safety restrictions~\cite{zhang2013ReachableDomainSpacecraft}, and ellipsoidal velocity-increment distributions~\cite{wen2018ReachableDomainSpacecraft}. This set-based view makes the connection between bounded maneuver authority and attainable spacecraft states explicit, providing a geometric description of where the spacecraft can go next.

The use of such information in control decisions also connects to a broader literature on constraint satisfaction. Viability theory examines the ability to maintain state constraints over time~\cite{aubin1991ViabilityTheory}, and verified feedback planning uses regions of attraction to organize admissible motion~\cite{tedrake2010LQRTrees}. In RL, safety has been addressed through both learning objectives and exploration strategies~\cite{garcia2015ComprehensiveSurveySafe}. Reachability-based approaches include iterative reachability estimation for policy optimization~\cite{ganai2023RESPO} and action projection using polynomial zonotopes~\cite{kochdumper2023ActionProjectionReachability}. The relationship between Hamilton--Jacobi reachability and RL is reviewed by Ganai et al.~\cite{ganai2024HJReachabilityRLSurvey}. Together, these perspectives motivate making control constraints part of the representation used for learned decisions.

Building on this viewpoint, this paper develops Reachability Analysis-Informed Reinforcement Learning (RARL) for deterministic multi-impulse interplanetary transfers. The central idea is to organize the learned decisions around intermediate waypoint selection, using local reachability geometry to define the available target set and deterministic astrodynamics to reconstruct the corresponding maneuvers. This places the selection of transfer geometry at the center of the policy's role while retaining dynamical reconstruction within an explicit physical model. The reachability description remains local, and the feasibility of the resulting nonlinear transfers is assessed empirically.

The contributions of this paper are fourfold. First, a reachability-informed action interface is formulated that relates bounded impulsive authority to node-dependent local waypoint sets. Second, a deterministic reconstruction layer connects policy-selected waypoints to nonlinear impulsive arcs, separating geometric decision making from maneuver reconstruction. Third, a terminal two-impulse proxy is formulated to assess the remaining rendezvous maneuver demand relative to the available control authority. Fourth, analytical and numerical studies examine the action mapping, the adequacy of the local reachability approximation, and the interpretation of the terminal proxy. These studies provide evidence for the principal modeling and interface choices in the proposed formulation.

The framework is evaluated on a deterministic Earth--Mars benchmark using the problem data of Marmo et al.~\cite{marmo2025ChanceConstraintRobustTrajectory}, with a separately computed sequential convex programming (SCP) solution serving as the numerical reference. Evaluation extends from nominal single-state training to a family of initial states induced by launch-time offsets, with the target state and transfer duration held fixed. Three independent training runs are performed for each regime, and the selected policies are assessed on a common held-out set of $10^4$ Monte Carlo (MC) departures. Reporting both individual policies and aggregate statistics distinguishes variation across independent training runs from variation across departure states within each policy. This establishes the consistency of the observed performance across the trained policies and quantifies the empirical robustness--cost tradeoff associated with extending policy reuse beyond the nominal departure condition.

The remainder of the paper is organized as follows. Section~\ref{sec:problem_formulation} defines the deterministic transfer problem, and Section~\ref{sec:methodology} develops the RARL formulation. Section~\ref{sec:validation} examines the action interface and its supporting approximations. Section~\ref{sec:case_study_results} presents the Earth--Mars configuration and the single-state and multi-state results. Section~\ref{sec:conclusions} concludes the paper.

%% file: sections/problem-formulation.tex
\section{Multi-Impulse Transfer Problem Formulation}
\label{sec:problem_formulation}

Consider a point-mass spacecraft with inertial position $\boldsymbol{r}(t)\in\mathbb{R}^{3}$ and velocity $\boldsymbol{v}(t)\in\mathbb{R}^{3}$. Its Cartesian state is
\begin{equation}
  \boldsymbol{x}(t)
  =
  \begin{bmatrix}
    \boldsymbol{r}(t) \\
    \boldsymbol{v}(t)
  \end{bmatrix}
  \in\mathbb{R}^{6}.
  \label{eq:state_definition}
\end{equation}
Between impulsive maneuvers, the spacecraft evolves according to the deterministic dynamics
\begin{equation}
  \dot{\boldsymbol{x}}(t)
  = f\!\left(\boldsymbol{x}(t)\right).
  \label{eq:continuous_dynamics}
\end{equation}
Here, $f(\cdot)$ denotes an autonomous deterministic force model. The Earth--Mars case study in Section~\ref{sec:case_study_results} uses a two-body heliocentric model consistently for propagation and maneuver reconstruction.

Let the fixed transfer interval $[t_0,t_f]$ be partitioned into $N$ ballistic arcs at the prescribed nodes $\{t_k\}_{k=0}^{N}$, where $t_N=t_f$ and $\Delta t_k=t_{k+1}-t_k$. The superscripts $-$ and $+$ denote the states immediately before and after an impulse, respectively. At the departure of arc $k$, an impulsive velocity change $\Delta\boldsymbol{v}_k\in\mathbb{R}^{3}$ is applied such that
\begin{equation}
  \boldsymbol{r}_k^{+}=\boldsymbol{r}_k^{-},
  \qquad
  \boldsymbol{v}_k^{+}
  =\boldsymbol{v}_k^{-}+\Delta\boldsymbol{v}_k,
  \qquad k=0,\ldots,N-1.
  \label{eq:impulse_jump}
\end{equation}
Denoting the deterministic flow map over an interval $\Delta t$ by $\mathcal{P}(\boldsymbol{x},\Delta t)$, each ballistic arc satisfies
\begin{equation}
  \boldsymbol{x}_{k+1}^{-}
  =
  \mathcal{P}\!\left(\boldsymbol{x}_k^{+},\Delta t_k\right),
  \qquad k=0,\ldots,N-1.
  \label{eq:ballistic_flow}
\end{equation}

An additional arrival impulse $\Delta\boldsymbol{v}_N$ is applied at the final node. It leaves the terminal position unchanged and closes the velocity rendezvous:
\begin{equation}
  \boldsymbol{r}_N^{+}=\boldsymbol{r}_N^{-},
  \qquad
  \boldsymbol{v}_N^{+}
  =\boldsymbol{v}_N^{-}+\Delta\boldsymbol{v}_N.
  \label{eq:arrival_impulse}
\end{equation}
For a prescribed initial state $\boldsymbol{x}_0$ and target state $\boldsymbol{x}_{T}(t_f)=[\boldsymbol{r}_{T}^{\mathsf{T}}(t_f), \boldsymbol{v}_{T}^{\mathsf{T}}(t_f)]^{\mathsf{T}}$, the boundary conditions are
\begin{equation}
  \boldsymbol{x}_0^{-}=\boldsymbol{x}_0,
  \qquad
  \boldsymbol{r}_N^{-}=\boldsymbol{r}_{T}(t_f),
  \qquad
  \boldsymbol{v}_N^{+}=\boldsymbol{v}_{T}(t_f).
  \label{eq:boundary_conditions}
\end{equation}

Each impulse is subject to the same magnitude limit,
\begin{equation}
  \left\|\Delta\boldsymbol{v}_k\right\|_2
  \leq \Delta V_{\max},
  \qquad k=0,\ldots,N,
  \label{eq:impulse_limit}
\end{equation}
and the impulsive transfer objective is
\begin{equation}
  \underset{\{\Delta\boldsymbol{v}_k\}_{k=0}^{N}}{\operatorname{minimize}}
  \quad
  J_{\Delta V}
  =
  \sum_{k=0}^{N}
  \left\|\Delta\boldsymbol{v}_k\right\|_2,
  \label{eq:impulsive_objective}
\end{equation}
subject to Eqs.~\eqref{eq:continuous_dynamics}--\eqref{eq:impulse_limit}. Thus, an $N$-arc transcription contains $N$ departure impulses and one terminal arrival impulse, giving $N+1$ impulse slots in both the objective and feasibility assessment.

%% file: sections/methodology.tex
\section{Reachability Analysis-Informed Reinforcement Learning Methodology}
\label{sec:methodology}

The proposed RARL framework combines learned waypoint selection with deterministic reconstruction of each transfer segment. At every policy-controlled node, a local state-transition map defines coordinates for the next position under bounded impulsive authority. The policy selects a waypoint in this coordinate system, after which a Lambert solve reconstructs the maneuver and ballistic arc required to reach it. A deterministic terminal solve then completes the position and velocity rendezvous.

\subsection{Moving-Target Representation and Local Sensitivity Map}
\label{subsec:moving_target_local_map}

Let $\boldsymbol{x}_{T}(t_f)$ denote the prescribed terminal target state. A time-aligned target state at node $t_k$ is obtained by propagating the terminal state backward over the remaining interval according to
\begin{equation}
  \boldsymbol{x}_{T}(t_k)
  =
  \mathcal{P}\!\left(
    \boldsymbol{x}_{T}(t_f),
    t_k-t_f
  \right),
  \qquad k=0,\ldots,N.
  \label{eq:moving_target}
\end{equation}
This moving reference enters the policy observation, while $\boldsymbol{x}_{T}(t_f)$ remains the prescribed terminal target state defined in Section~\ref{sec:problem_formulation}.

The local action interface follows from the sensitivity of a one-segment coast to a velocity perturbation. Consider the nominal coast from the current pre-impulse state,
\begin{equation}
  \boldsymbol{x}_{k+1}^{\mathrm{nom}}
  =
  \mathcal{P}\!\left(
    \boldsymbol{x}_{k}^{-},
    \Delta t_k
  \right).
  \label{eq:nominal_flow}
\end{equation}
The state transition matrix (STM) over the same interval is
\begin{equation}
  \boldsymbol{\Phi}_k
  =
  \left.
  \frac{\partial \mathcal{P}(\boldsymbol{x},\Delta t_k)}
       {\partial \boldsymbol{x}}
  \right|_{\boldsymbol{x}=\boldsymbol{x}_k^{-}}
  =
  \begin{bmatrix}
    \boldsymbol{\Phi}_{rr,k} & \boldsymbol{\Phi}_{rv,k} \\
    \boldsymbol{\Phi}_{vr,k} & \boldsymbol{\Phi}_{vv,k}
  \end{bmatrix}.
  \label{eq:stm_partition}
\end{equation}
For a purely impulsive perturbation at node $k$,
\begin{equation}
  \delta\boldsymbol{x}_k
  =
  \begin{bmatrix}
    \boldsymbol{0} \\
    \delta\boldsymbol{v}_k
  \end{bmatrix},
  \qquad
  \delta\boldsymbol{x}_{k+1}
  \approx
  \boldsymbol{\Phi}_k\delta\boldsymbol{x}_k.
  \label{eq:impulsive_state_perturbation}
\end{equation}
The first-order position deviation at the next node therefore satisfies
\begin{equation}
  \Delta\boldsymbol{r}_{k+1}
  =
  \boldsymbol{r}_{k+1}-\boldsymbol{r}_{k+1}^{\mathrm{nom}}
  \approx
  \boldsymbol{\Phi}_{rv,k}\,\delta\boldsymbol{v}_k.
  \label{eq:local_position_map}
\end{equation}
Here, $\delta\boldsymbol{v}_k$ denotes a local perturbation used to construct the sensitivity map and the notation $\Delta\boldsymbol{v}_k$ is reserved for a reconstructed physical maneuver. The matrix $\boldsymbol{A}_k=\boldsymbol{\Phi}_{rv,k}\in\mathbb{R}^{3\times3}$ is the local position-from-velocity sensitivity map. Its singular vectors and singular values characterize the directions and amplification with which local velocity perturbations map into next-node position deviations.

\subsection{Reachable-Waypoint Action Interface}
\label{subsec:waypoint_interface}

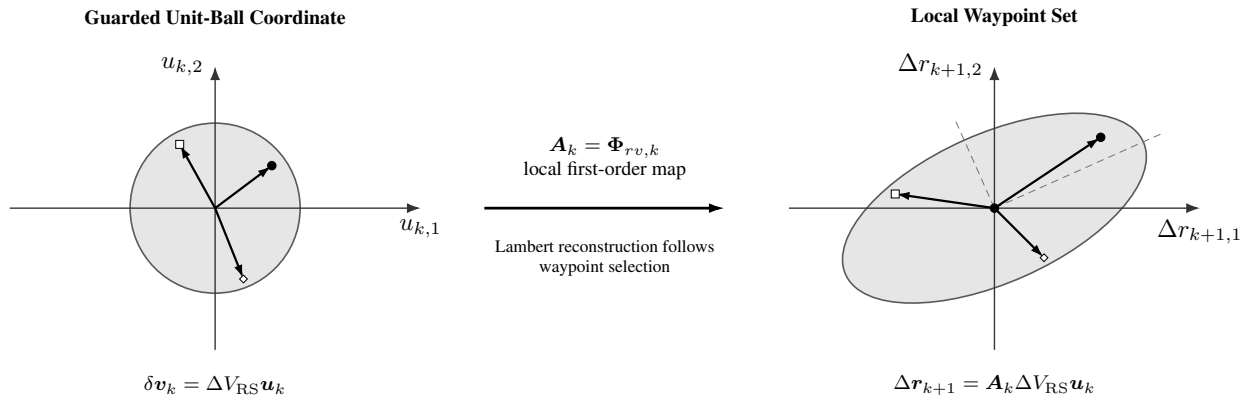
\begin{figure}[htbp]
  \centering
  \resizebox{\linewidth}{!}{%
    \input{figs/local-reachability-mapping.tex}%
  }
  \caption{Two-dimensional schematic of the local first-order map. The position-from-velocity STM block maps guarded unit-ball coordinates to an ellipsoidal waypoint set centered on the nominal coast position, with orientation and semiaxis lengths set by local segment geometry.}
  \label{fig:local_reachability_map}
\end{figure}

The sensitivity relation defines a waypoint interface once the local velocity perturbation is bounded. Let the closed Euclidean unit ball be $\mathbb{B}_2=\{\boldsymbol{u}\in\mathbb{R}^3\mid\lVert\boldsymbol{u}\rVert_2\leq1\}$. The physical impulse limit remains $\Delta V_{\max}$, while the linear waypoint interface uses the contracted radius
\begin{equation}
  \Delta V_{\mathrm{RS}}
  =
  \Delta V_{\max}-g,
  \qquad
  0\leq g<\Delta V_{\max},
  \label{eq:guarded_interface_radius}
\end{equation}
where $g$ is an inner guard. The corresponding first-order next-position set is
\begin{equation}
  \mathcal{R}_{k}^{\mathrm{lin}}
  =
  \left\{
    \boldsymbol{r}_{k+1}^{\mathrm{nom}}
    +
    \boldsymbol{A}_{k}
    \left(\Delta V_{\mathrm{RS}}\boldsymbol{u}\right)
    \middle|
    \boldsymbol{u}\in\mathbb{B}_2
  \right\}.
  \label{eq:linear_reachable_set}
\end{equation}
The set $\mathcal{R}_{k}^{\mathrm{lin}}$ is the affine image of $\mathbb{B}_2$ under $\boldsymbol{A}_k$ and therefore forms an ellipsoid centered at $\boldsymbol{r}_{k+1}^{\mathrm{nom}}$. If $\boldsymbol{A}_k$ is rank deficient, this image becomes a degenerate ellipsoid~\cite{kurzhanski2000ellipsoidal,wen2018ReachableDomainSpacecraft}. Its principal directions and semiaxis lengths vary with the current orbital geometry, segment duration and adopted dynamics through $\boldsymbol{A}_k$. Figure~\ref{fig:local_reachability_map} illustrates how the guarded unit-ball coordinate is transformed by $\boldsymbol{A}_k$ into this locally defined waypoint geometry.

\pagebreak

The policy produces $\boldsymbol{a}_k\in[-1,1]^3$. After componentwise clipping to this box, the action is mapped radially into $\mathbb{B}_2$ according to
\begin{equation}
  \boldsymbol{u}_k
  =
  \begin{cases}
    \boldsymbol{0},
      & \lVert\boldsymbol{a}_k\rVert_2\leq\varepsilon, \\[3pt]
    \left(\min\{1,\lVert\boldsymbol{a}_k\rVert_2\}\right)^p
    \dfrac{\boldsymbol{a}_k}{\lVert\boldsymbol{a}_k\rVert_2},
      & \text{otherwise},
  \end{cases}
  \label{eq:radial_action_map}
\end{equation}
where $\varepsilon>0$ is a numerical zero threshold and $p>0$ is a radial exponent. Values $p>1$ shift the represented radius toward the interior, whereas values $p<1$ shift it toward the boundary. All experiments reported in this paper use $p=1$. A zero policy action selects the nominal coast position, while actions outside the unit ball are projected radially onto its boundary.

The selected waypoint is
\begin{equation}
  \boldsymbol{r}_{k+1}^{\mathrm{RS}}
  =
  \boldsymbol{r}_{k+1}^{\mathrm{nom}}
  +
  \boldsymbol{A}_k
  \left(\Delta V_{\mathrm{RS}}\boldsymbol{u}_k\right).
  \label{eq:reachable_waypoint}
\end{equation}
For a fixed spacecraft state and segment duration, Eq.~\eqref{eq:reachable_waypoint} maps $\boldsymbol{u}_k$ affinely to the selected waypoint, whereas the preceding radial transformation from the box-constrained policy action $\boldsymbol{a}_k$ to $\boldsymbol{u}_k$ is generally nonlinear.

The guard provides an empirical margin for small nonlinear reconstruction discrepancies observed at the action boundary, and its selection and validation are reported in Sections~\ref{sec:validation} and~\ref{sec:case_study_results}. However, it does not guarantee that every Lambert-reconstructed impulse satisfies the physical limit $\Delta V_{\max}$. The first-order set is therefore used as a local decision interface rather than as a nonlinear reachability or recursive-feasibility certificate.

\subsection{Deterministic Segment Reconstruction and Terminal Closure}
\label{subsec:segment_reconstruction}

For $k=0,\ldots,N-2$, the policy-selected waypoint is realized through a zero-revolution prograde Lambert solve~\cite{izzo2015RevisitingLambert} over $\Delta t_k$,
\begin{equation}
  \left(
    \boldsymbol{v}_{k,\mathrm{dep}}^{\mathrm{pro}},
    \boldsymbol{v}_{k+1,\mathrm{arr}}^{\mathrm{pro}}
  \right)
  =
  \operatorname{Lambert}_{0,\mathrm{pro}}\!\left(
    \boldsymbol{r}_k,
    \boldsymbol{r}_{k+1}^{\mathrm{RS}},
    \Delta t_k
  \right).
  \label{eq:waypoint_lambert}
\end{equation}
The resulting departure impulse and next pre-impulse state are
\begin{subequations}
\label{eq:waypoint_reconstruction}
\begin{align}
  \Delta\boldsymbol{v}_k
  & =
  \boldsymbol{v}_{k,\mathrm{dep}}^{\mathrm{pro}}
  -\boldsymbol{v}_k^{-},
  \label{eq:waypoint_reconstructed_impulse} \\
  \boldsymbol{x}_{k+1}^{-}
  & =
  \begin{bmatrix}
    \boldsymbol{r}_{k+1}^{\mathrm{RS}} \\
    \boldsymbol{v}_{k+1,\mathrm{arr}}^{\mathrm{pro}}
  \end{bmatrix}.
  \label{eq:waypoint_state_update}
\end{align}
\end{subequations}
This reconstruction makes the selected waypoint dynamically consistent with the adopted two-body model when the Lambert solution exists. Compliance with the physical impulse limit is evaluated separately because membership in the linear set does not guarantee that the Lambert-reconstructed impulse obeys the same limit.

The policy controls the first $N-1$ transfers. At node $N-1$, a zero-revolution Lambert rendezvous closes both terminal position and velocity. The prograde and retrograde branches are evaluated. Let $\mathcal{B}_{N-1}\subseteq\{\mathrm{pro},\mathrm{retro}\}$ contain the branches for which the Lambert solver returns a solution. When this set is nonempty, the branch with the lower total terminal impulse is selected. For $b\in\mathcal{B}_{N-1}$,
\begin{subequations}
\label{eq:terminal_branch_selection}
\begin{align}
  \left(
    \boldsymbol{v}_{N-1,\mathrm{dep}}^{0,b},
    \boldsymbol{v}_{N,\mathrm{arr}}^{0,b}
  \right)
  & =
  \operatorname{Lambert}_{0,b}\!\left(
    \boldsymbol{r}_{N-1},\!
    \boldsymbol{r}_{T}(t_f),\!
    \Delta t_{N-1}
  \right),
  \label{eq:terminal_branch_velocities} \\
  b^{\star}
  & =
  \underset{b\in\mathcal{B}_{N-1}}{\operatorname{arg\,min}}
  \left(
    \left\lVert\Delta\boldsymbol{v}_{N-1}^{(b)}\right\rVert_2
    \right.
    \notag \\
  & \quad
    \left.
    +
    \left\lVert\Delta\boldsymbol{v}_{N}^{(b)}\right\rVert_2
  \right),
  \label{eq:terminal_branch_choice} \\
  \Delta\boldsymbol{v}_{N-1}^{(b)}
  & =
  \boldsymbol{v}_{N-1,\mathrm{dep}}^{0,b}
  -\boldsymbol{v}_{N-1}^{-},
  \label{eq:terminal_departure_impulse} \\
  \Delta\boldsymbol{v}_{N}^{(b)}
  & =
  \boldsymbol{v}_{T}(t_f)
  -\boldsymbol{v}_{N,\mathrm{arr}}^{0,b}.
  \label{eq:terminal_arrival_impulse_method}
\end{align}
\end{subequations}
If $\mathcal{B}_{N-1}$ is empty, the transition is truncated as a terminal Lambert failure. Accordingly, the reconstructed trajectory contains $N-1$ intermediate departure impulses and two terminal impulses, consistent with the $N+1$ impulse slots defined in Section~\ref{sec:problem_formulation}. When the selected terminal solve succeeds, the environment assigns $\boldsymbol{x}_N=\boldsymbol{x}_T(t_f)$, so the terminal error is zero by construction under the adopted dynamics model rather than through a property learned by the policy or a guarantee under unmodeled dynamics.

\pagebreak

\subsection{MDP Observation and Reward Construction}
\label{subsec:mdp_reward}

RARL is formulated as a finite-horizon Markov decision process (MDP). At each node, the observation augments the absolute spacecraft state with its position and velocity relative to the moving target in Eq.~\eqref{eq:moving_target}. These relative quantities are expressed in the spacecraft-defined radial-transverse-normal (RTN) frame. The radial direction is aligned with $\boldsymbol{r}_k$, the normal direction is aligned with $\boldsymbol{r}_k\times\boldsymbol{v}_k$, and the transverse direction completes the right-handed orthonormal triad. Let $\boldsymbol{C}_{\mathcal{I}\rightarrow\mathrm{RTN}}(\boldsymbol{x}_k^{-})$ denote the corresponding inertial-to-RTN rotation. The target-relative position and velocity differences expressed in this RTN frame are
\begin{subequations}
\label{eq:relative_state_rtn}
\begin{align}
  \Delta\boldsymbol{r}_k^{\mathrm{RTN}}
  & =
  \boldsymbol{C}_{\mathcal{I}\rightarrow\mathrm{RTN}}(
    \boldsymbol{x}_k^{-})
  \left(
    \boldsymbol{r}_{T}(t_k)-\boldsymbol{r}_k^{-}
  \right),
  \label{eq:relative_position_rtn} \\
  \Delta\boldsymbol{v}_k^{\mathrm{RTN}}
  & =
  \boldsymbol{C}_{\mathcal{I}\rightarrow\mathrm{RTN}}(
    \boldsymbol{x}_k^{-})
  \left(
    \boldsymbol{v}_{T}(t_k)-\boldsymbol{v}_k^{-}
  \right).
  \label{eq:relative_velocity_rtn}
\end{align}
\end{subequations}

The $14$-dimensional observation is
\begin{equation}
  \boldsymbol{o}_k
  =
  \mathcal{S}\!\left(
  \begin{bmatrix}
    \boldsymbol{x}_k^{-} \\
    \Delta\boldsymbol{r}_k^{\mathrm{RTN}} \\
    \Delta\boldsymbol{v}_k^{\mathrm{RTN}} \\
    \tau_k \\
    \rho_{k-1}
  \end{bmatrix}
  \right),
  \qquad
  \tau_k=\frac{t_f-t_k}{t_f-t_0},
  \label{eq:observation_definition}
\end{equation}
where $\mathcal{S}(\cdot)$ applies fixed componentwise affine scaling and clips the resulting components to $[-1,1]$. The previous impulse ratio is
\begin{equation}
  \rho_{k-1}
  =
  \frac{
    \left\lVert\Delta\boldsymbol{v}_{k-1}\right\rVert_2
  }{\Delta V_{\max}},
  \qquad
  \rho_{-1}=0.
  \label{eq:previous_impulse_ratio}
\end{equation}
Section~\ref{sec:case_study_results} defines the characteristic units and fixed scaling bounds used in the present implementation.

A low-cost two-impulse assessment can be defined at a generic node over the remaining time $\Delta t_{k\rightarrow f}=t_f-t_k$. The nominal terminal state and its mismatch from the prescribed target are
\begin{subequations}
\label{eq:remaining_terminal_mismatch}
\begin{align}
  \boldsymbol{x}_{f}^{\mathrm{nom}}
  & =
  \mathcal{P}\!\left(
    \boldsymbol{x}_k^{-},
    \Delta t_{k\rightarrow f}
  \right),
  \label{eq:remaining_nominal_terminal_state} \\
  \Delta\boldsymbol{x}_{k\rightarrow f}
  & =
  \boldsymbol{x}_{T}(t_f)-\boldsymbol{x}_{f}^{\mathrm{nom}}
  =
  \begin{bmatrix}
    \Delta\boldsymbol{r}_{k\rightarrow f} \\
    \Delta\boldsymbol{v}_{k\rightarrow f}
  \end{bmatrix}.
  \label{eq:remaining_terminal_state_error}
\end{align}
\end{subequations}
Over the same remaining interval, let $\boldsymbol{\Phi}_{k\rightarrow f}$ denote the STM, partitioned as in Eq.~\eqref{eq:stm_partition}. Using the Moore--Penrose pseudoinverse~\cite{horn2012matrix,golub2013matrix}, denoted by $(\cdot)^{\dagger}$, the minimum-norm departure estimate and the corresponding arrival correction are
\begin{subequations}
\label{eq:two_impulse_estimates}
\begin{align}
  \Delta\boldsymbol{v}_{1,k}^{\mathrm{lin}}
  & =
  \boldsymbol{\Phi}_{rv,k\rightarrow f}^{\dagger}
  \Delta\boldsymbol{r}_{k\rightarrow f},
  \label{eq:two_impulse_departure_estimate} \\
  \Delta\boldsymbol{v}_{2,k}^{\mathrm{lin}}
  & =
  \Delta\boldsymbol{v}_{k\rightarrow f}
  -
  \boldsymbol{\Phi}_{vv,k\rightarrow f}
  \Delta\boldsymbol{v}_{1,k}^{\mathrm{lin}}.
  \label{eq:two_impulse_arrival_estimate}
\end{align}
\end{subequations}
Although $\boldsymbol{\Phi}_{rv,k\rightarrow f}$ is square, it may become singular or ill-conditioned for particular orbital geometries and times-to-go. The pseudoinverse, evaluated using the prescribed singular-value cutoff, therefore provides a numerically stable minimum-norm least-squares construction. When the matrix is nonsingular and all its singular values are retained, the pseudoinverse reduces to the ordinary inverse. Since the underlying sensitivity relation is first order, the resulting corrections remain local estimates of the required maneuver authority. Their magnitudes are normalized by $\Delta V_{\max}$, and the associated one-sided excess is defined by
\begin{subequations}
\label{eq:two_impulse_violation}
\begin{align}
  \eta_{1,k}
  & =
  \frac{\left\lVert
    \Delta\boldsymbol{v}_{1,k}^{\mathrm{lin}}
  \right\rVert_2}{\Delta V_{\max}},
  \label{eq:two_impulse_departure_ratio} \\
  \eta_{2,k}
  & =
  \frac{\left\lVert
    \Delta\boldsymbol{v}_{2,k}^{\mathrm{lin}}
  \right\rVert_2}{\Delta V_{\max}},
  \label{eq:two_impulse_arrival_ratio} \\
  \xi_k
  & =
  \left[\max\{\eta_{1,k},\eta_{2,k}\}-1\right]_{+},
  \label{eq:two_impulse_excess}
\end{align}
\end{subequations}
where $[z]_{+}=\max\{0,z\}$. The construction can also be evaluated at earlier nodes when an anticipatory learning signal is desired. A value $\xi_k>0$ at such a node does not establish infeasibility of the multi-impulse problem because additional intermediate maneuvers remain available. The present implementation evaluates $\xi_k$ only at $k=N-1$, where the remaining horizon consists of the final segment. At this node, $\xi_{N-1}$ is zero when both estimated corrections remain within $\Delta V_{\max}$ and increases only when either estimate exceeds that limit. The resulting one-sided measure informs learning only when the estimated terminal maneuver demand exceeds the available per-impulse authority. Within that authority, candidate terminal solutions remain differentiated by their reconstructed Lambert impulse costs.

An alternative full-state least-squares construction is presented in Appendix~\ref{app:secondary_proxy}. It estimates a single linearized correction from the combined terminal position and velocity mismatch. The two-impulse measure is adopted in the present formulation because its departure-and-arrival structure mirrors the terminal rendezvous reconstruction. The full-state formulation is retained as a complementary diagnostic that may be useful when a single-correction surrogate is preferred.

For a unified reward definition, let $\Delta\boldsymbol{v}_{k,1}$ denote the departure impulse on transition $k$ and let $\Delta\boldsymbol{v}_{k,2}$ denote the additional arrival impulse on the terminal transition. Define the per-impulse cap exceedance as
\begin{equation}
  \epsilon_{k,i}^{\Delta v}
  =
  \left[
    \left\lVert\Delta\boldsymbol{v}_{k,i}\right\rVert_2
    -\Delta V_{\max}
  \right]_{+},
  \qquad i\in\{1,2\}.
  \label{eq:impulse_exceedance}
\end{equation}
For a named event $\mathcal{E}$, let $\indicator{\mathcal{E}}$ equal one when the event occurs and zero otherwise. The terminal-transition indicator $\indicator{\{k=N-1\}}$ is one only for the final segment. Let $\mathcal{E}_k^{\mathrm{L}}$ denote failure of the Lambert reconstruction on transition $k$, whether on an intermediate or terminal segment, and let $\mathcal{E}_k^{\mathrm{U}}$ denote detection of an unbound two-body state on that transition. The four reward contributions account for reconstructed maneuver expenditure, per-impulse cap exceedance, final-segment proxy shaping and truncation-triggering failures, respectively. The transition reward is decomposed as
\begin{subequations}
\label{eq:rarl_reward}
\begin{align}
  r_k
  & =
  r_k^{\Delta v}
  +r_k^{\mathrm{cap}}
  +r_k^{\mathrm{proxy}}
  +r_k^{\mathrm{fail}},
  \label{eq:rarl_reward_total} \\
  r_k^{\Delta v}
  & =
  -\lambda_{\Delta v}
  \left(
    \left\lVert\Delta\boldsymbol{v}_{k,1}\right\rVert_2
    +
    \indicator{\{k=N-1\}}
    \left\lVert\Delta\boldsymbol{v}_{k,2}\right\rVert_2
  \right),
  \label{eq:rarl_reward_maneuver} \\
  r_k^{\mathrm{cap}}
  & =
  -\lambda_{\mathrm{viol}}
  \left(
    \epsilon_{k,1}^{\Delta v}
    +
    \indicator{\{k=N-1\}}
    \epsilon_{k,2}^{\Delta v}
  \right),
  \label{eq:rarl_reward_cap} \\
  r_k^{\mathrm{proxy}}
  & =
  -\indicator{\{k=N-1\}}
  \lambda_{\eta}\xi_{N-1},
  \label{eq:rarl_reward_proxy} \\
  r_k^{\mathrm{fail}}
  & =
  -\lambda_{\mathrm{L}}
  \indicator{\mathcal{E}_k^{\mathrm{L}}}
  -
  \lambda_{\mathrm{U}}
  \indicator{\mathcal{E}_k^{\mathrm{U}}}.
  \label{eq:rarl_failure_reward}
\end{align}
\end{subequations}
The nonnegative weights $\lambda_{\Delta v}$, $\lambda_{\mathrm{viol}}$, $\lambda_{\eta}$, $\lambda_{\mathrm{L}}$ and $\lambda_{\mathrm{U}}$ scale the maneuver, cap-exceedance, proxy, Lambert-failure and unbound-state contributions, respectively. The proxy contribution is active only on the final transition. It vanishes when both linearized correction estimates remain within $\Delta V_{\max}$ and increases in magnitude according to the larger normalized exceedance. Unlike the cap-exceedance contribution, this term depends on the STM-based two-impulse assessment rather than the impulses returned by the terminal Lambert reconstruction.

An episode terminates normally after $N$ segments and is truncated if an intermediate or terminal Lambert solve fails or if an unbound two-body state is detected. An impulse-limit exceedance contributes the cap term in Eq.~\eqref{eq:rarl_reward} but does not itself truncate the episode. The impulse costs and cap exceedances are evaluated in physical velocity units. The dimensional cap exceedance is related to the impulse-to-limit ratio through
\begin{equation}
  \epsilon_{k,i}^{\Delta v}
  =
  \Delta V_{\max}
  \left[
    \frac{\left\lVert\Delta\boldsymbol{v}_{k,i}\right\rVert_2}
         {\Delta V_{\max}}
    -1
  \right]_{+}.
  \label{eq:dimensional_ratio_equivalence}
\end{equation}
Thus, dimensional and ratio-based penalties impose the same threshold after the associated weight is rescaled. The dimensional form in Eq.~\eqref{eq:impulse_exceedance} is retained so that the maneuver and cap-exceedance terms are both expressed in physical velocity units. This decomposition places maneuver expenditure, impulse-limit exceedance, final-segment proxy shaping and truncating numerical or dynamical failures within a single transition reward.

\subsection{PPO Training Workflow}
\label{subsec:ppo_workflow}

The policy $\pi_{\theta}(\boldsymbol{a}_k\mid\boldsymbol{o}_k)$ is represented by a stochastic actor--critic model and trained using proximal policy optimization (PPO)~\cite{schulman2017ProximalPolicyOptimization}. The actor neural network (NN) parameterizes this conditional action distribution, while the critic NN approximates the expected discounted return from $\boldsymbol{o}_k$ through the value function $V_{\phi}(\boldsymbol{o}_k)$. Section~\ref{sec:case_study_results} defines the network dimensions, PPO hyperparameters, rollout budget and checkpoint-selection rules used in the present implementation.

Algorithm~\ref{alg:rarl_workflow} summarizes the interaction between policy learning and deterministic reconstruction. Here, $N_{\mathrm{roll}}$ denotes the prescribed number of transitions collected in each fixed-length rollout before a PPO update. The generic workflow may stop when the prescribed training budget is exhausted or when a prespecified convergence criterion is satisfied. The reported experiments use the fixed training budgets specified in Section~\ref{sec:case_study_results}.

\begin{algorithm}[htbp]
  \caption{RARL Rollout and Training Workflow}
  \label{alg:rarl_workflow}
  \begin{algorithmic}[1]
    \Require \parbox[t]{0.91\linewidth}{$\boldsymbol{x}_0$, $\boldsymbol{x}_T(t_f)$, $\{t_k\}_{k=0}^{N}$, $\Delta V_{\max}$, $g$, $p$, $\mathcal{P}(\cdot,\Delta t)$, Lambert solver, actor and critic NN architectures, reward weights, $N_{\mathrm{roll}}$ and a prespecified stopping criterion}
    \Ensure Trained RARL policy $\pi_\theta^{\star}$
    \State \textbf{Initialization:} Initialize the actor parameters $\theta$,
    critic parameters $\phi$ and PPO optimizer
    \Repeat
      \State \textbf{Rollout collection:} Initialize an empty rollout buffer
      $\mathcal{D}$
      \While{$\lvert\mathcal{D}\rvert<N_{\mathrm{roll}}$}
        \State \textbf{Episode initialization:} If required, initialize the
        episode at $\boldsymbol{x}_0$ with $k\gets0$ and $\rho_{-1}\gets0$
        \State \textbf{Observation:} Evaluate the moving target using
        Eq.~\eqref{eq:moving_target} and form $\boldsymbol{o}_k$ using
        Eq.~\eqref{eq:observation_definition}
        \If{$k\leq N-2$}
          \State \textbf{Policy query:} Sample $\boldsymbol{a}_k$ from
          $\pi_\theta$
          \State \textbf{Reachability map:} Propagate the nominal coast
          and STM
          \State \textbf{Action projection:} Map $\boldsymbol{a}_k$ using
          Eq.~\eqref{eq:radial_action_map}
          \State \textbf{Reachable-waypoint selection:}
          Construct
          $\boldsymbol{r}_{k+1}^{\mathrm{RS}}$ using
          Eq.~\eqref{eq:reachable_waypoint}
          \State \textbf{Lambert reconstruction:} Solve the prograde Lambert
          reconstruction
          \State \textbf{State update:} If the Lambert reconstruction succeeds,
          update $\boldsymbol{x}_{k+1}^{-}$ and $\rho_k$. Otherwise, set
          $\mathcal{E}_k^{\mathrm{L}}$
        \Else
          \State \textbf{Terminal closure:} Propagate the terminal nominal
          coast and form the terminal STM
          \State Evaluate $\xi_{N-1}$ using Eq.~\eqref{eq:two_impulse_violation}
          \State Solve the terminal branches using Eq.~\eqref{eq:terminal_branch_selection}
          \State If $\mathcal{B}_{N-1}$ is empty, set
          $\mathcal{E}_{N-1}^{\mathrm{L}}$. Otherwise, assign
          $\boldsymbol{x}_N\gets\boldsymbol{x}_T(t_f)$
        \EndIf
        \State \textbf{Failure assessment:} Set $\mathcal{E}_k^{\mathrm{U}}$
        if the state is unbound
        \State \textbf{Reward and transition:} Evaluate Eq.~\eqref{eq:rarl_reward}
        and store the transition in $\mathcal{D}$
        \State \textbf{Episode status:} Truncate on
        $\mathcal{E}_k^{\mathrm{L}}$ or $\mathcal{E}_k^{\mathrm{U}}$, terminate
        if $k=N-1$, and otherwise set $k\gets k+1$
      \EndWhile
      \State \textbf{Policy update:} Update $\theta$ and $\phi$ with PPO using
      $\mathcal{D}$
    \Until{budget exhausted or convergence criterion met}
    \State \Return $\pi_\theta^{\star}$
  \end{algorithmic}
\end{algorithm}
\FloatBarrier

%% file: figs/local-reachability-mapping.tex
\begin{tikzpicture}[
  >=Latex,
  axis/.style={->,draw=black!80,line width=0.55pt},
  set/.style={draw=black!65,fill=black!10,line width=0.65pt},
  vector/.style={-{Latex[length=2.0mm,width=1.2mm]},draw=black,line width=0.9pt},
  map/.style={-{Latex[length=2.3mm,width=1.4mm]},draw=black,line width=1.1pt},
  note/.style={fill=white,inner sep=1.5pt,font=\footnotesize},
  sample/.style={circle,draw=black,fill=black,inner sep=1.2pt}
]
  \begin{scope}[shift={(-55mm,0)}]
    \node[font=\footnotesize\bfseries] at (0,27mm)
      {Guarded Unit-Ball Coordinate};
    \draw[set] (0,0) circle (12mm);
    \draw[axis] (-29mm,0) -- (29mm,0) node[below=1pt] {$u_{k,1}$};
    \draw[axis] (0,-20mm) -- (0,20mm) node[left=1pt] {$u_{k,2}$};
    \coordinate (u1) at (8mm,6mm);
    \coordinate (u2) at (-5mm,9mm);
    \coordinate (u3) at (4mm,-10mm);
    \draw[vector] (0,0) -- (u1);
    \draw[vector] (0,0) -- (u2);
    \draw[vector] (0,0) -- (u3);
    \node[sample] at (u1) {};
    \node[rectangle,draw=black,fill=white,minimum size=3.2pt,inner sep=0pt] at (u2) {};
    \node[diamond,draw=black,fill=white,minimum size=3.3pt,inner sep=0pt] at (u3) {};
    \node[note] at (0,-25mm)
      {$\delta\boldsymbol{v}_{k}=\Delta V_{\mathrm{RS}}\boldsymbol{u}_k$};
  \end{scope}

  \begin{scope}[shift={(55mm,0)}]
    \node[font=\footnotesize\bfseries] at (0,27mm)
      {Local Waypoint Set};
    \draw[set,rotate=24] (0,0) ellipse [x radius=23mm,y radius=10.5mm];
    \draw[axis] (-29mm,0) -- (29mm,0) node[below=1pt] {$\Delta r_{k+1,1}$};
    \draw[axis] (0,-20mm) -- (0,20mm) node[left=1pt] {$\Delta r_{k+1,2}$};
    \draw[black!55,densely dashed] (0,0) -- (24:26mm);
    \draw[black!55,densely dashed] (0,0) -- (114:14mm);
    \coordinate (y1) at (15mm,10mm);
    \coordinate (y2) at (-14mm,2mm);
    \coordinate (y3) at (7mm,-7mm);
    \draw[vector] (0,0) -- (y1);
    \draw[vector] (0,0) -- (y2);
    \draw[vector] (0,0) -- (y3);
    \node[sample] at (y1) {};
    \node[rectangle,draw=black,fill=white,minimum size=3.2pt,inner sep=0pt] at (y2) {};
    \node[diamond,draw=black,fill=white,minimum size=3.3pt,inner sep=0pt] at (y3) {};
    \node[sample] at (0,0) {};
    \node[note] at (0,-25mm)
      {$\Delta\boldsymbol{r}_{k+1}=\boldsymbol{A}_k\Delta V_{\mathrm{RS}}\boldsymbol{u}_k$};
  \end{scope}

  \draw[map] (-17mm,0) -- (17mm,0);
  \node[note,align=center] at (0,7mm)
    {$\boldsymbol{A}_k=\boldsymbol{\Phi}_{rv,k}$\\local first-order map};
  \node[font=\scriptsize,align=center] at (0,-7mm)
    {Lambert reconstruction follows\\waypoint selection};
\end{tikzpicture}

%% file: sections/validation.tex
\section{Analysis and Validation of the Reachability Interface}
\label{sec:validation}

This section analyzes the action mapping, evaluates the accuracy of the local reachability map, and examines the terminal two-impulse proxy as a measure of maneuver demand. The studies use the deterministic Earth--Mars configuration described in Section~\ref{sec:case_study_results}.

\subsection{Action-Mapping Properties}
\label{subsec:validation_action_mapping}

The cube-to-ball transformation determines how policy actions populate the local waypoint set. Boundedness ensures that the selected point belongs to $\mathcal{R}_k^{\mathrm{lin}}$, while centering, directional behavior, and radial scaling determine how coast, interior, and boundary selections are represented. The centered $L_2$ map in Eq.~\eqref{eq:radial_action_map}, with $p=1$, is compared with two alternative parameterizations under identical uniform cube inputs.

The first alternative is an equal-volume spherical construction. Let $q_i=(a_i+1)/2$ for $i=1,2,3$. The mapped coordinate is
\begin{equation}
  \boldsymbol{u}
  =
  \begin{bmatrix}
    s\sqrt{1-z^2}\cos\phi \\
    s\sqrt{1-z^2}\sin\phi \\
    sz
  \end{bmatrix},
  \quad
  \begin{aligned}
    s &= q_1^{1/3}, \\
    z &= 1-2q_2, \\
    \phi &= 2\pi q_3.
  \end{aligned}
  \label{eq:validation_equal_volume_map}
\end{equation}
Uniform cube inputs produce a uniform distribution in the unit ball, but $\boldsymbol{a}=\boldsymbol{0}$ maps to $\boldsymbol{u}=[-2^{-1/3},0,0]^{\mathsf{T}}$ rather than the coast center. The second alternative is the $L_\infty$-based radial map,
\begin{equation}
  \boldsymbol{u}
  =
  \frac{\boldsymbol{a}}{\lVert\boldsymbol{a}\rVert_2}
  \lVert\boldsymbol{a}\rVert_\infty,
  \qquad \lVert\boldsymbol{a}\rVert_2>\varepsilon,
  \label{eq:validation_linf_map}
\end{equation}
with $\boldsymbol{u}=\boldsymbol{0}$ otherwise. Like the centered $L_2$ map, it preserves zero-action centering and the direction of nonzero inputs. Both radial maps therefore retain the angular distribution induced by the cube, while differing in the radius assigned to each input.

\begin{figure}[htbp]
  \centering
  \includegraphics[width=5.5in,trim=0 6bp 0 6bp,clip]{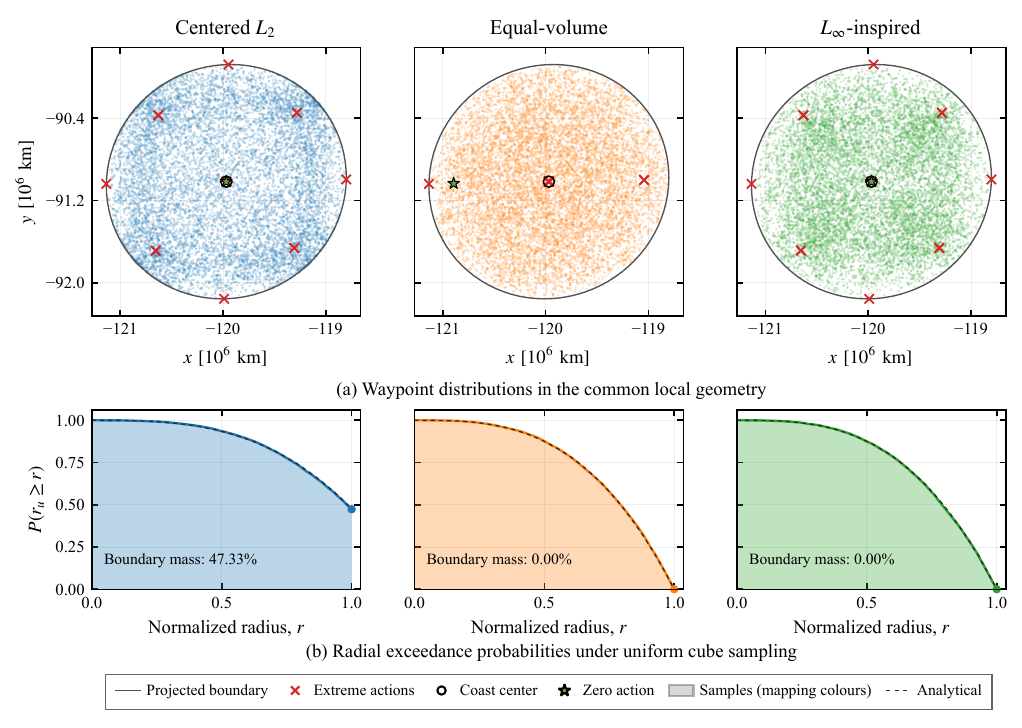}
  \caption{Action-mapping comparison for $10^4$ uniform cube inputs.
  (a) Projected waypoints at the initial Earth--Mars node.
  (b) Radial exceedance probabilities for $r_u=\|\boldsymbol{u}\|_2$.
  Dashed curves show analytical values, and annotations report empirical boundary probabilities.}
  \label{fig:action_mapping_diagnostics}
\end{figure}

Figure~\ref{fig:action_mapping_diagnostics} compares $10^4$ uniformly sampled cube inputs passed through each map. The projected waypoints use the same STM evaluated at the initial benchmark state and the same guarded radius. These controlled inputs isolate the effect of the parameterization from changes in the local orbital geometry.

For the ideal maps, with the numerical zero threshold set to zero, the inclusive radial exceedance probabilities for $0\leq r\leq1$ are
\begin{equation}
  \Pr(r_u\geq r)
  =
  \begin{cases}
    1-\pi r^3/6, & \text{centered }L_2, \\[2pt]
    1-r^3, & \text{either alternative}.
  \end{cases}
  \label{eq:validation_radial_exceedance}
\end{equation}
As the centered $L_2$ map leaves inputs inside the unit ball unchanged and projects the remaining cube volume onto its boundary, the probability inside radius $r$ follows from the ball-to-cube volume ratio. The theoretical boundary probability is therefore $\Pr(r_u=1)=1-\pi/6\approx47.64\%$, closely matching the empirical value of $47.33\%$, as stated in Fig.~\ref{fig:action_mapping_diagnostics}. Neither alternative mapping produced boundary samples in this comparison. The empirical and analytical exceedance curves distinguish the radial behavior of the centered $L_2$ and $L_\infty$ mappings, despite their similar projected waypoint distributions.

In their ideal forms, both radial maps are continuous and differentiable almost everywhere. The centered map changes derivative at the unit sphere, while the $L_\infty$ map also changes derivative where the dominant action component changes. The centered map combines a direct local relation $\boldsymbol{u}=\boldsymbol{a}$ with a finite input volume assigned to the boundary. These properties motivate its use for representing coast, interior, and high-authority waypoint selections. Boundary selection refers to the linear waypoint set, with the physical impulse subsequently determined by Lambert reconstruction.

\subsection{Accuracy of the Local Linear Map}
\label{subsec:validation_linearization}

The first-order map is assessed by comparing its predicted position change with nonlinear propagation from the same initial state. For a local impulse perturbation $\delta\boldsymbol{v}_k$, define
\begin{equation}
  \boldsymbol{x}_{k+1}^{\mathrm{true}}(\delta\boldsymbol{v}_k)
  =
  \mathcal{P}\!\left(
    \begin{bmatrix}
      \boldsymbol{r}_k^{-} \\
      \boldsymbol{v}_k^{-}+\delta\boldsymbol{v}_k
    \end{bmatrix},\Delta t_k
  \right).
  \label{eq:validation_true_next_state}
\end{equation}
Its position component gives the nonlinear displacement from the nominal coast in Eq.~\eqref{eq:nominal_flow}. The nonlinear and linearized displacements are
\begin{subequations}
\begin{align}
  \Delta\boldsymbol{r}_{k+1}^{\mathrm{true}}
  &=\boldsymbol{r}_{k+1}^{\mathrm{true}}
    -\boldsymbol{r}_{k+1}^{\mathrm{nom}},
  \label{eq:validation_true_dr} \\
  \Delta\boldsymbol{r}_{k+1}^{\mathrm{lin}}
  &=\boldsymbol{A}_k\delta\boldsymbol{v}_k.
  \label{eq:validation_lin_dr}
\end{align}
\end{subequations}
The absolute and relative errors are then
\begin{subequations}
\begin{align}
  e_r
  &=\left\lVert
    \Delta\boldsymbol{r}_{k+1}^{\mathrm{true}}
    -\Delta\boldsymbol{r}_{k+1}^{\mathrm{lin}}
  \right\rVert_2,
  \label{eq:validation_lin_abs_error} \\
  \bar e_r
  &=\frac{e_r}
    {\left\lVert\Delta\boldsymbol{r}_{k+1}^{\mathrm{true}}\right\rVert_2}.
  \label{eq:validation_lin_rel_error}
\end{align}
\end{subequations}
The relative error is normalized by the magnitude of the impulse-induced position change.

\begin{figure}[htbp]
  \centering
  \captionsetup{width=4.510in}
  \includegraphics[width=3.630in]{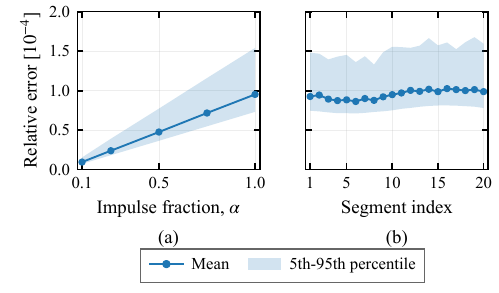}
  \caption{Linearization adequacy of the first-order reachability map.
  (a) Relative error versus impulse fraction.
  (b) Segment-wise error for $\alpha=1$ (maximum impulse).}
  \label{fig:linearization_adequacy}
\end{figure}

The diagnostic uses 20 consecutive nominal coast nodes from the benchmark initial state, with a segment duration of $17.4395$ days. At each node, 200 independent directions are sampled for each of five impulse fractions,
$\alpha\in\{0.1,0.25,0.5,0.75,1.0\}$, with
\begin{equation}
  \delta\boldsymbol{v}_k=\alpha\Delta V_{\max}\boldsymbol{d},
  \qquad
  \boldsymbol{d}\sim\mathcal{U}(\mathbb{S}^{2}),
  \label{eq:validation_sampled_impulse}
\end{equation}
where $\mathbb{S}^{2}$ is the unit sphere. Figure~\ref{fig:linearization_adequacy} shows the mean relative error and the 5th--95th-percentile range across sampled directions and nodes. At the physical cap, the mean relative error is $9.52\times10^{-5}$, the 95th percentile is $1.54\times10^{-4}$, and the maximum is $1.81\times10^{-4}$. These results support the local position approximation over the tested coast states and impulse amplitudes in the considered Earth--Mars case study.

The selected waypoints are realized through nonlinear Lambert reconstruction as described in Section~\ref{subsec:segment_reconstruction}. Checks of the reconstructed intermediate impulses showed no cap violations over the confirmation runs. The inner guard is treated as an empirical margin and does not guarantee containment.

\subsection{Analysis of the Terminal Two-Impulse Proxy}
\label{subsec:validation_terminal_proxy}

\begin{figure}[htbp]
  \centering
  \captionsetup{width=4.510in}
  \includegraphics[width=3.300in]{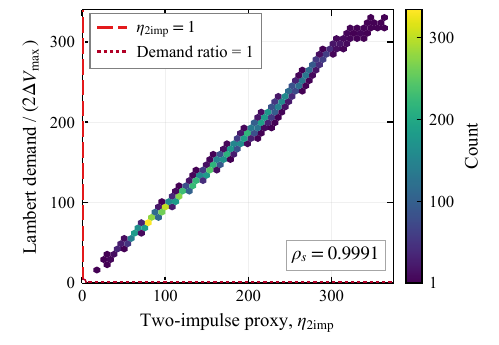}
  \caption{Terminal two-impulse proxy versus normalized Lambert demand. Hexbins aggregate $10^4$ infeasible samples. Dashed and dotted lines mark activation and total-demand thresholds, respectively.}
  \label{fig:eta2imp_validation}
\end{figure}

The terminal two-impulse proxy assesses maneuver demand using the first-order departure and arrival corrections defined in Section~\ref{sec:methodology}. At the penultimate node, let $\eta_{2\mathrm{imp}}=\max\{\eta_{1,N-1},\eta_{2,N-1}\}$ denote the larger normalized estimate. The reward contribution in Eq.~\eqref{eq:rarl_reward_proxy} then depends on $\xi_{N-1}=[\eta_{2\mathrm{imp}}-1]_+$. Below the activation threshold, terminal solutions remain differentiated by their reconstructed impulse costs. Above it, the assessment provides a graded measure of estimated excess demand. Figure~\ref{fig:eta2imp_validation} examines how this proxy ranks nonlinear terminal maneuver demand.

The terminal diagnostic comprises $10^4$ penultimate states generated by random impulsive coasts from the benchmark initial condition. For this diagnostic, the two-impulse estimates in Eq.~\eqref{eq:two_impulse_estimates} were computed using the Yamanaka--Ankersen STM~\cite{yamanaka2002NewStateTransition}. The nonlinear terminal demand was obtained by solving the zero-revolution prograde and retrograde Lambert branches and selecting the successful solution with minimum total departure-and-arrival impulse cost. The plotted demand is this total normalized by the combined allowance $2\Delta V_{\max}$ for the departure and arrival impulses. A normalized total at or below unity satisfies the combined impulse allowance, although either maneuver may still exceed its individual cap. Terminal feasibility therefore also depends on how the demand is distributed between the departure and arrival impulses.

The Spearman rank correlation~\cite{spearman1904proof,zar1999biostatistical} is $\rho_s=0.9991$, indicating a strong monotone association within the sampled set. All samples are terminal-infeasible and have $\eta_{2\mathrm{imp}}>1$, placing the comparison in the regime where the one-sided penalty is active. Since $\xi_{N-1}=\eta_{2\mathrm{imp}}-1$ throughout this set, the excess used in the reward has exactly the same ranking. The reference lines identify the proxy activation threshold and the combined impulse budget. The observed association supports two-impulse assessment as a measure of terminal demand within the sampled infeasible regime.

%% file: sections/case-study-results.tex
\section{Earth--Mars Transfer Case Study}
\label{sec:case_study_results}

\subsection{Benchmark and Numerical Configuration}
\label{subsec:cs_problem_setup}

The case study considers a deterministic Earth--Mars rendezvous based on the benchmark of Marmo et al.~\cite{marmo2025ChanceConstraintRobustTrajectory}. The prescribed initial and target states are expressed in the heliocentric ecliptic J2000 frame. Between maneuver nodes, the spacecraft evolves under two-body dynamics,
\begin{equation}
  \dot{\boldsymbol{r}}=\boldsymbol{v},
  \qquad
  \dot{\boldsymbol{v}}
  =-\mu_{\odot}\frac{\boldsymbol{r}}{\lVert\boldsymbol{r}\rVert_2^3},
  \label{eq:em_two_body_dynamics}
\end{equation}
where $\mu_{\odot}$ is the solar gravitational parameter. This force model is used consistently for propagation, STM construction, and Lambert reconstruction. Two-body propagation and Lambert reconstruction in the RARL environment use PyKEP~\cite{izzo2026pykep}.

\begin{table}[htbp]
  \centering
  \caption{Earth--Mars benchmark parameters and numerical scales.}
  \label{tab:em_problem_data}
  \footnotesize
  \begin{tabular}{@{}l c l@{}}
    \toprule
    Quantity & Value & Unit \\
    \midrule
    $\mu_{\odot}$ & $132{,}712{,}440{,}018$ & $\mathrm{km^3/s^2}$ \\
    $t_f-t_0$ & $348.79$ & days \\
    $N$ & $20$ & -- \\
    $\Delta t_k$ & $17.4395$ & days \\
    $\Delta V_{\max}$ & $0.7600$ & $\mathrm{km/s}$ \\
    \midrule
    $\boldsymbol{r}_0^{-}$ &
    $[-140{,}699{,}693,\,-51{,}614{,}428,\,980]^{\mathsf{T}}$ & km \\
    $\boldsymbol{v}_0^{-}$ &
    $[9.7746,\,-28.0783,\,4.3377\times10^{-4}]^{\mathsf{T}}$ & $\mathrm{km/s}$ \\
    $\boldsymbol{r}_T(t_f)$ &
    $[-172{,}682{,}023,\,176{,}959{,}469,\,7{,}948{,}912]^{\mathsf{T}}$ & km \\
    $\boldsymbol{v}_T(t_f)$ &
    $[-16.4274,\,-14.8605,\,9.2149\times10^{-2}]^{\mathsf{T}}$ & $\mathrm{km/s}$ \\
    \midrule
    $r_{\mathrm{scale}}$ & $1.4960\times10^8$ & km \\
    $v_{\mathrm{scale}}$ & $29.7845$ & $\mathrm{km/s}$ \\
    $t_{\mathrm{scale}}$ & $5.02275\times10^6$ & s \\
    \bottomrule
  \end{tabular}
\end{table}

Table~\ref{tab:em_problem_data} lists the nominal boundary states and numerical configuration. The fixed transfer duration $t_f-t_0$ is divided into $N=20$ ballistic segments of equal duration $\Delta t_k$, giving 21 impulse slots including the terminal arrival impulse. The initial position $\boldsymbol{r}_0^{-}$ and velocity $\boldsymbol{v}_0^{-}$ are specified before the departure impulse. The target position $\boldsymbol{r}_T(t_f)$ and velocity $\boldsymbol{v}_T(t_f)$ specify the required state after the arrival impulse, consistent with Eq.~\eqref{eq:boundary_conditions}. Terminal rendezvous is imposed through the reconstruction described in Section~\ref{subsec:segment_reconstruction}.

The environment uses nondimensional position, velocity, and time, with characteristic scales $r_{\mathrm{scale}}$, $v_{\mathrm{scale}}=\sqrt{\mu_{\odot}/r_{\mathrm{scale}}}$, and $t_{\mathrm{scale}}=r_{\mathrm{scale}}/v_{\mathrm{scale}}$. These choices give a unit nondimensional gravitational parameter. Reported distances, maneuver magnitudes, and times are expressed in physical units. The inner guard introduced in Section~\ref{subsec:waypoint_interface} is set to $g=2\times10^{-4}\,\mathrm{km/s}$ in Eq.~\eqref{eq:guarded_interface_radius}.

\subsection{Training and Evaluation Protocol}
\label{subsec:cs_training_protocol}

RARL policies are trained either at the nominal departure state, termed single-state training, or across a prescribed family of initial states, termed multi-state training. Both regimes use the 14-dimensional observation and three-dimensional action defined in Section~\ref{sec:methodology}. Positions and velocities are first nondimensionalized using the scales in Table~\ref{tab:em_problem_data}. Table~\ref{tab:em_obs_bounds} lists the fixed componentwise bounds used by $\mathcal{S}(\cdot)$ in Eq.~\eqref{eq:observation_definition} throughout training and evaluation. The reward weights in Eq.~\eqref{eq:rarl_reward} are $\lambda_{\Delta v}=25$, $\lambda_{\mathrm{viol}}=200$, $\lambda_{\eta}=10$, and $\lambda_{\mathrm{L}}=\lambda_{\mathrm{U}}=500$ for both training regimes. Online reward normalization is applied during training, while evaluation uses the unnormalized rewards and physical maneuver costs.

\begin{table}[htbp]
  \centering
  \caption{Bounds for fixed observation scaling.}
  \label{tab:em_obs_bounds}
  \footnotesize
  \begin{tabular}{@{}l c c@{}}
    \toprule
    Input & Lower bound & Upper bound \\
    \midrule
    $\boldsymbol{r}_k^{-}/r_{\mathrm{scale}}$ &
    $[-2.3,-2.4,-0.3]$ & $[2.3,2.4,0.3]$ \\
    $\boldsymbol{v}_k^{-}/v_{\mathrm{scale}}$ &
    $[-2.1,-2.2,-0.25]$ & $[2.1,2.2,0.25]$ \\
    $\Delta\boldsymbol{r}_k^{\mathrm{RTN}}/r_{\mathrm{scale}}$ &
    $[-3.5,-2.7,-0.4]$ & $[3.5,2.7,0.4]$ \\
    $\Delta\boldsymbol{v}_k^{\mathrm{RTN}}/v_{\mathrm{scale}}$ &
    $[-1.6,-2.5,-0.2]$ & $[1.6,2.5,0.2]$ \\
    $\tau_k$ & $0$ & $1$ \\
    $\rho_{k-1}$ & $0$ & $2$ \\
    \bottomrule
  \end{tabular}
\end{table}

The actor and critic use separate multilayer perceptrons (MLPs), each with three hidden layers and hyperbolic tangent activations. The hidden layer sizes are listed in Table~\ref{tab:em_ppo_config}. Together with the output layers, the model contains approximately $5{,}600$ trainable parameters. PPO uses generalized advantage estimation (GAE), and the learning rate and clip range decrease linearly over the prescribed training budget. Single-state and multi-state training use the same network and PPO settings, with three independently initialized policies per regime. Each single-state run uses approximately 240 million training transitions, with evaluation every $9{,}600$ transitions. Each multi-state run uses approximately 160 million transitions, with evaluation every $96{,}000$ transitions. Stable-Baselines3~\cite{stable-baselines3} provides the PPO implementation, built on PyTorch~\cite{paszke2019pytorch,ansel2024pytorch2}, while the training environments are implemented as custom classes compatible with Gymnasium~\cite{towers2025gymnasium}.

\begin{table}[htbp]
  \centering
  \caption{Policy architecture and PPO training configuration.}
  \label{tab:em_ppo_config}
  \footnotesize
  \begin{tabular}{@{}l c c@{}}
    \toprule
    Parameter & Symbol & Value \\
    \midrule
    Actor hidden sizes & -- & $[70,32,15]$ \\
    Critic hidden sizes & -- & $[56,14,4]$ \\
    Initial action log-std & $\log\sigma_{\pi,0}$ & $0$ \\
    Initial learning rate & $\alpha_0$ & $3.0\times10^{-4}$ \\
    Final learning rate & $\alpha_{\mathrm{final}}$ & $1.0\times10^{-5}$ \\
    Discount factor & $\gamma$ & $0.9999$ \\
    GAE factor & $\lambda_{\mathrm{GAE}}$ & $0.99$ \\
    Initial clip range & $\varepsilon_{\mathrm{clip},0}$ & $0.20$ \\
    Final clip range & $\varepsilon_{\mathrm{clip,final}}$ & $0.10$ \\
    Entropy coefficient & $c_{\mathrm{ent}}$ & $0.01$ \\
    Value-function coefficient & $c_V$ & $0.6$ \\
    Maximum gradient norm & -- & $0.5$ \\
    Target Kullback--Leibler divergence & -- & $0.05$ \\
    Vectorized environments & $N_{\mathrm{env}}$ & $8$ \\
    Rollout transitions per environment & $n_{\mathrm{steps}}$ & $1{,}200$ \\
    Total rollout transitions & $N_{\mathrm{roll}}$ & $9{,}600$ \\
    Minibatch size & $B_{\mathrm{mb}}$ & $150$ \\
    PPO epochs per update & $n_{\mathrm{epochs}}$ & $10$ \\
    \bottomrule
  \end{tabular}
\end{table}

Single-state policies are trained and evaluated at the nominal departure state. Multi-state training uses a launch-window-inspired initial-state family obtained by propagating the nominal departure state through offsets $\Delta t_{\mathrm{dep}}\in[-2,2]$ days under the same two-body model. The terminal target state and the transfer duration remain fixed. Table~\ref{tab:em_state_sets} summarizes the three state sets. The training set contains 100 states on a uniform offset grid within the prescribed interval, and each training episode samples one of these states uniformly. The separate ten-state evaluation grid is interleaved between training states and has no overlap with the training grid. It assesses interpolation and is used for policy checkpoint selection. The MC validation set contains $10^4$ states with offsets sampled uniformly from $[-2,2]$ days. It is used only after policy selection to estimate maneuver-cost distributions and feasibility rates over the prescribed family. Each state set is shared across the three policies.

\begin{table}[htbp]
  \centering
  \caption{State sets for training, evaluation, and MC validation.}
  \label{tab:em_state_sets}
  \footnotesize
  \begin{tabular}{@{}>{\raggedright\arraybackslash}p{0.27\columnwidth} c >{\raggedright\arraybackslash}p{0.49\columnwidth}@{}}
    \toprule
    Set & States & Construction and purpose \\
    \midrule
    Training & $100$ & Uniform grid sampled during training \\
    Interleaved \mbox{evaluation} & $10$ & Separate grid between training states for checkpoint selection \\
    MC validation & $10^4$ & Uniformly sampled offsets for assessment after selection \\
    \bottomrule
  \end{tabular}
\end{table}

During evaluation, actions are selected deterministically from the mean of the policy distribution and passed through the action mapping in Eq.~\eqref{eq:radial_action_map}. An episode is considered feasible if it completes the transfer without a Lambert failure or an unbound-state truncation and all 21 impulses satisfy the physical cap, allowing a numerical comparison tolerance of $10^{-12}\,\mathrm{km/s}$. Policy checkpoints are ranked first by the fraction of feasible episodes, then by the fraction without truncation, the mean maximum impulse-cap violation, and the mean total maneuver cost. The first two quantities are maximized and the latter two minimized, with violations within the numerical tolerance treated as zero. Among fully feasible candidates, this selects the lowest mean-cost checkpoint. Selection uses the nominal state for single-state training and the ten-state interleaved evaluation set for multi-state training. The selected policies are retained for the reported results.

Results are reported for each trained policy and summarized by the mean and sample standard deviation across three independent training runs with different random seeds. For the multi-state results, variation across MC validation states is reported separately from variation across independently trained policies.

\subsection{Deterministic SCP Reference}
\label{subsec:cs_scp_reference}

A deterministic reference trajectory is computed using SCP, which solves a sequence of convex approximations to the nonlinear transfer problem~\cite{malyuta2022convex}. The calculation uses the boundary states, transfer duration, discretization and physical impulse cap in Table~\ref{tab:em_problem_data}, with the total maneuver objective defined in Eq.~\eqref{eq:impulsive_objective}. All 21 impulses, including the terminal arrival impulse, are optimization variables.

The state and impulse variables are nondimensionalized using the scales in Section~\ref{subsec:cs_problem_setup}. At iteration $j$, the impulse application and subsequent two-body propagation over segment $k$ are approximated by
\begin{equation}
  \widetilde{\boldsymbol{x}}_{k+1}^{-}
  =\boldsymbol{C}_{k}^{(j)}\widetilde{\boldsymbol{x}}_{k}^{-}
  +\boldsymbol{D}_{k}^{(j)}\Delta\widetilde{\boldsymbol{v}}_{k}
  +\boldsymbol{d}_{k}^{(j)},
  \label{eq:scp_affine_dynamics}
\end{equation}
for $k=0,\ldots,N-1$. Here, tildes denote nondimensional quantities, $\boldsymbol{C}_{k}^{(j)}\in\mathbb{R}^{6\times6}$ and $\boldsymbol{D}_{k}^{(j)}\in\mathbb{R}^{6\times3}$ are the state and impulse Jacobians evaluated along the reference trajectory, and $\boldsymbol{d}_{k}^{(j)}$ ensures agreement with nonlinear propagation at that reference. The state Jacobian is obtained by integrating the two-body variational equations, and the impulse Jacobian comprises its three velocity columns. The state and STM are integrated together using the DOP853 solver in SciPy~\cite{virtanen2020scipy}, with relative and absolute tolerances of $10^{-12}$ and $10^{-13}$, respectively.

Each convex subproblem minimizes the sum of impulse magnitudes subject to the affine dynamics, the prescribed initial and terminal conditions, and the impulse bounds. Euclidean trust regions constrain changes from the reference state and impulse at each node, with nondimensional radii of $0.1$ and $0.02$, respectively. A numerical margin of $10^{-9}\,\mathrm{km/s}$ is subtracted from the physical impulse cap within the optimizer. The resulting second-order cone programs are formulated in CVXPY~\cite{diamond2016cvxpy,agrawal2018rewriting} and solved using MOSEK~\cite{mosek11}, with primal feasibility, dual feasibility and relative-gap tolerances of $10^{-10}$.

The initialization originates from a two-impulse Lambert transfer between the prescribed endpoint positions, with zero intermediate impulses. A stored SCP solution obtained from this initialization is subsequently refined to the tolerances specified here. After each convex solve, the updated impulses are propagated from the prescribed initial state under the nonlinear dynamics to construct the next reference. Convergence requires terminal position and velocity errors no greater than $10^{-3}\,\mathrm{km}$ and $10^{-6}\,\mathrm{km/s}$, respectively, with no physical impulse-cap exceedance. It also requires a maximum absolute componentwise difference of $10^{-11}$ between the optimized and nonlinearly propagated nondimensional node states, and a change in total maneuver cost no greater than $10^{-8}\,\mathrm{km/s}$.

The nominal refinement converges in four SCP iterations. Independent propagation of the resulting impulse sequence using PyKEP~\cite{izzo2026pykep} gives terminal position and velocity errors of $2.19\times10^{-4}\,\mathrm{km}$ and $2.20\times10^{-11}\,\mathrm{km/s}$, respectively, with all impulses within the physical cap. This solution provides the validated local SCP reference used in the subsequent performance comparison.

\subsection{Single-State Benchmark Performance}
\label{subsec:results_single_state}

\begin{table}[htbp]
\centering
\begin{threeparttable}
\caption{Nominal single-state performance and impulse allocation.}
\label{tab:em_single_state_results}
\small
\setlength{\tabcolsep}{3pt}
\begin{tabular*}{\textwidth}{@{\extracolsep{\fill}}llccccc@{}}
\toprule
Quantity & Unit & SCP & \multicolumn{4}{c}{RARL (single-state)} \\
\cmidrule(lr){4-7}
 & & & Seed 1 & Seed 2 & Seed 3 & Mean $\pm$ SD\tnote{a} \\
\midrule
\multicolumn{7}{@{}l}{\textit{Transfer performance}} \\
Total maneuver cost & km/s & $10.0585$ & $10.2069$ & $10.2469$ & $10.2416$ & $10.2318 \pm 0.0218$ \\
Cost above SCP\tnote{b} & km/s & $0.0000$ & $0.1484$ & $0.1884$ & $0.1831$ & $0.1733 \pm 0.0218$ \\
Gap to SCP\tnote{b} & \% & $0.00$ & $1.47$ & $1.87$ & $1.82$ & $1.72 \pm 0.22$ \\
\midrule
\multicolumn{7}{@{}l}{\textit{Impulse allocation}} \\
Near-coast impulses\tnote{c} & count & 7 & 5 & 5 & 6 & $5.3 \pm 0.6$ \\
Near-coast maneuver cost\tnote{c} & km/s & $0.0000$ & $0.0627$ & $0.0702$ & $0.1347$ & $0.0892 \pm 0.0396$ \\
Near-cap impulses\tnote{d} & count & 13 & 12 & 12 & 12 & 12 \\
Terminal departure impulse & km/s & $0.7600$ & $0.7424$ & $0.7248$ & $0.7354$ & $0.7342 \pm 0.0089$ \\
Terminal arrival impulse & km/s & $0.7600$ & $0.7541$ & $0.7427$ & $0.7411$ & $0.7460 \pm 0.0071$ \\
\bottomrule
\end{tabular*}
\begin{tablenotes}[flushleft]
\small
\setlength{\itemsep}{0pt}
\setlength{\parsep}{0pt}
\item[a] Mean $\pm$ sample standard deviation (SD) over three independently trained policies. Identical counts are shown once.
\item[b] Absolute and percentage cost differences are relative to the SCP reference.
\item[c] Near-coast counts impulses with $\|\Delta\boldsymbol{v}_k\|_2<0.05\,\mathrm{km/s}$.
\item[d] Near-cap counts impulses with $\|\Delta\boldsymbol{v}_k\|_2>0.95\Delta V_{\max}$.
\end{tablenotes}
\end{threeparttable}
\end{table}

The single-state results correspond to the nominal departure condition, $\Delta t_{\mathrm{dep}}=0$. Table~\ref{tab:em_single_state_results} compares the three selected RARL policies with the SCP reference described in Section~\ref{subsec:cs_scp_reference}. All three RARL policies complete the nominal transfer without truncation or impulse-cap violations. Terminal rendezvous is imposed through the Lambert reconstruction described in Section~\ref{subsec:segment_reconstruction}.

\begin{figure}[htbp]
  \centering
  \captionsetup{width=4.510in}
  \includegraphics[width=3.500in]{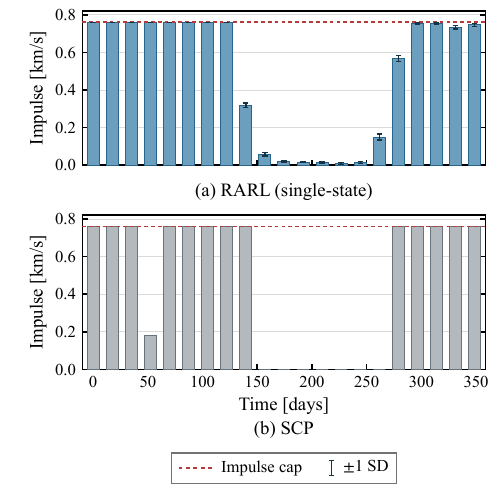}
  \caption{Nominal impulse allocation for (a) RARL and (b) SCP. RARL bars and whiskers show the mean and sample standard deviation across three independent runs.}
  \label{fig:res_ctrl_single}
\end{figure}

\pagebreak

The total maneuver cost is $10.2318\pm0.0218\,\mathrm{km/s}$, reported as the mean $\pm$ sample standard deviation across the three independently trained policies. This corresponds to a mean cost difference of $0.1733\,\mathrm{km/s}$, or $1.72\%$, above the SCP reference of $10.0585\,\mathrm{km/s}$. Marmo et al.~\cite{marmo2025ChanceConstraintRobustTrajectory} report a deterministic optimal control solution with a total maneuver cost of $10.4048\,\mathrm{km/s}$, placing the present RARL costs close to the reported literature result. These results show that the proposed RARL methodology recovers benchmark-quality nominal performance across three independently trained policies, with costs close to the SCP reference. Presenting each policy alongside the aggregate statistics makes the observed variation across independent training runs explicit, both in total maneuver cost and in impulse allocation. The consistently feasible nominal outcomes and small cost dispersion provide evidence that the reported performance is repeatable across the tested training seeds.

Figure~\ref{fig:res_ctrl_single} compares the node-wise RARL allocation with the SCP solution. The RARL bars and whiskers report the mean and sample standard deviation at each impulse node. Both profiles exhibit an approximate bang--off--bang structure, with large early and late maneuvers separated by a near-coast portion. This broad allocation pattern is also evident in the reported deterministic optimal control solution~\cite{marmo2025ChanceConstraintRobustTrajectory}. Using the thresholds defined in Table~\ref{tab:em_single_state_results}, each RARL policy has 12 near-cap impulses and 5 or 6 near-coast impulses, compared with 13 and 7, respectively, for SCP. The RARL near-coast impulses contribute $0.0892\pm0.0396\,\mathrm{km/s}$ across the three policies. The terminal departure and arrival impulses remain below the physical cap, with respective mean values of $0.7342\pm0.0089\,\mathrm{km/s}$ and $0.7460\pm0.0071\,\mathrm{km/s}$.

\begin{figure}[htbp]
  \centering
  \includegraphics[width=\linewidth,trim=0 0 0 14bp,clip]{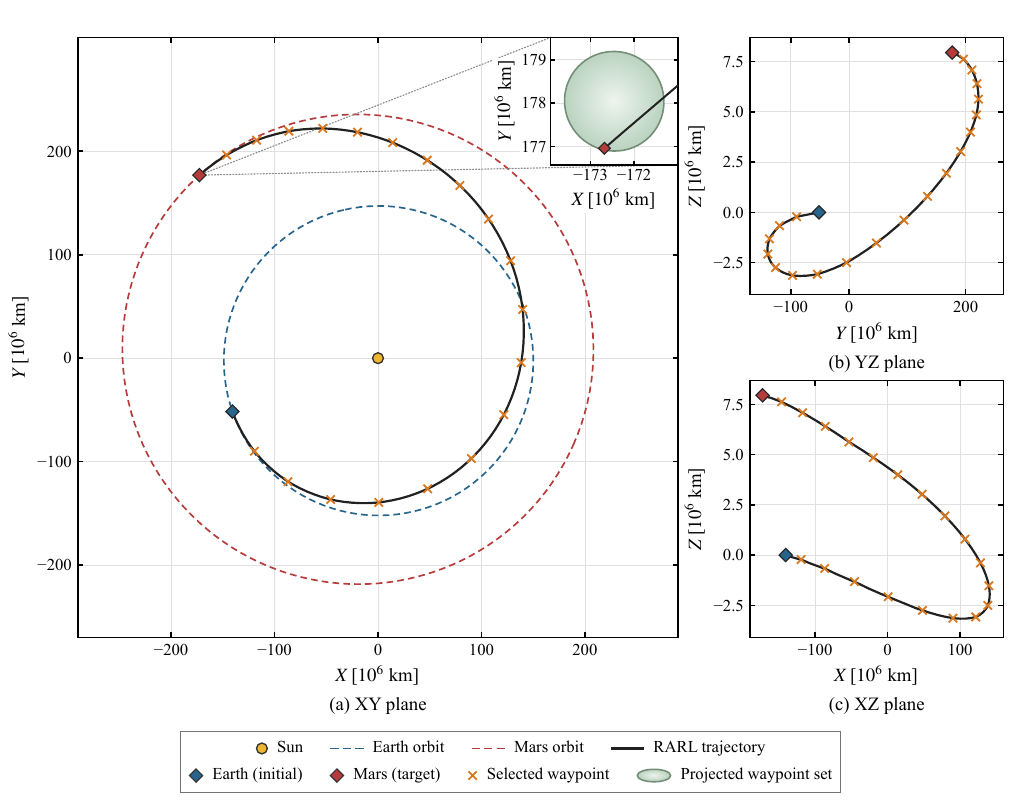}
  \caption{Nominal Earth--Mars transfer for the representative single-state RARL policy. Panel (a) shows the $XY$ projection, with a terminal inset showing the projected waypoint set. Panels (b)--(c) show the $YZ$ and $XZ$ projections. Crosses mark selected intermediate waypoints.}
  \label{fig:res_traj_single}
\end{figure}

Figure~\ref{fig:res_traj_single} shows a representative single-state RARL trajectory, generated by the Seed 3 policy, whose nominal maneuver cost of $10.2416\,\mathrm{km/s}$ is the median of the three runs, as reported in Table~\ref{tab:em_single_state_results}. The three spatial projections show the transfer geometry, while the terminal inset resolves the final approach to the prescribed target.

\begin{figure}[H]
  \centering
  \includegraphics[width=\linewidth,height=0.92\textheight,keepaspectratio,trim=0 0 0 6bp,clip]{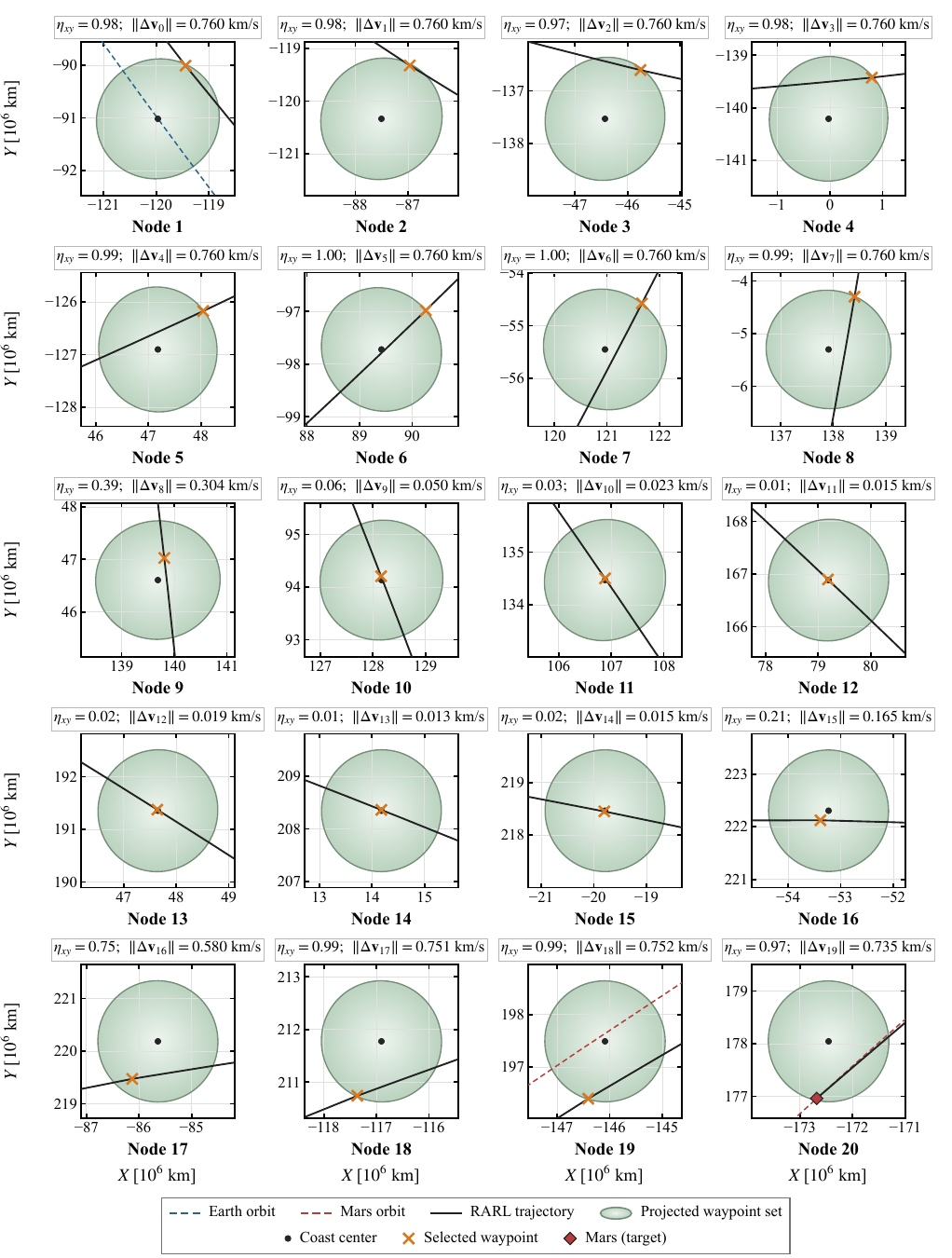}
  \caption{Local $XY$ projections of the waypoint sets along the representative nominal trajectory. Crosses mark policy-selected waypoints at nodes 1--19, and the diamond marks the terminal target at node 20. Annotations give the normalized projected radius $\eta_{xy}$ and the preceding reconstructed impulse magnitude. All panels use the same spatial scale.}
  \label{fig:res_rs_single}
\end{figure}

Figure~\ref{fig:res_rs_single} presents the corresponding local $XY$ projections of the waypoint sets at every arrival node along the trajectory. Nodes 1--19 show policy-selected waypoints, and node 20 shows the prescribed terminal target. Each panel pairs the arrival-node geometry with the impulse applied at the preceding departure node.

The normalized projected radius $\eta_{xy}$ measures the displacement from the coast center relative to the ellipse boundary in the same projected direction. It is zero at the center and unity on the boundary. The accompanying impulse magnitude is obtained from the three-dimensional Lambert reconstruction. Together, these quantities describe the waypoint's location within the projected set and the maneuver required to realize it.

At nodes 1--8, the selected waypoints lie near the projected boundary, with $\eta_{xy}$ between approximately $0.97$ and $1.00$, and the preceding impulses are close to the physical cap. Node 9 marks the transition toward the interior, with $\eta_{xy}=0.386$ and an impulse of $0.304\,\mathrm{km/s}$. At nodes 10--15, $\eta_{xy}$ lies between $0.014$ and $0.064$, while the corresponding impulses range from approximately $0.013$ to $0.050\,\mathrm{km/s}$ and contribute $0.1347\,\mathrm{km/s}$ in total. The selections then move outward through nodes 16--17 and return close to the boundary at nodes 18--19. This progression connects the near-coast portion and subsequent increase in maneuver effort to the local waypoint geometry exposed to the policy.

The nominal results therefore demonstrate that waypoint selection coupled with deterministic reconstruction can achieve maneuver costs close to the SCP reference while preserving a direct geometric interpretation of the policy decisions.

\subsection{Multi-State Performance under Launch-Time Dispersion}
\label{subsec:results_multi_state}

Section~\ref{subsec:results_single_state} establishes the nominal-transfer performance of RARL. Multi-state evaluation examines whether this capability extends to dispersed departure conditions, allowing a trained policy to be reused across transfers with different initial states. Launch-time offsets in $[-2,2]$ days produce dispersion in the initial position and velocity, while the target state and transfer duration remain those of the nominal benchmark. The policies must therefore accommodate the resulting changes in transfer geometry through waypoint selection, seeking low maneuver cost while satisfying terminal rendezvous and the individual impulse limits.

The state sets in Table~\ref{tab:em_state_sets} separate policy training, selection on interleaved evaluation states, and subsequent MC validation. The selected single-state and multi-state policies are assessed on the same $10^4$ MC departures, providing a common basis for comparing feasibility and maneuver-cost distributions. Each selected policy is applied across these departures without further training. Table~\ref{tab:em_multi_state_results} compares the selected single-state and multi-state RARL policies over this departure family and at the nominal departure state.

\begin{table}[htbp]
\centering
\begin{threeparttable}
\caption{Multi-state performance under launch-time dispersion. All quantities except feasibility rate and percentage gaps are in km/s.}
\label{tab:em_multi_state_results}
\small
\setlength{\tabcolsep}{3pt}
\begin{tabular*}{\textwidth}{@{\extracolsep{\fill}}lccccc@{}}
\toprule
Quantity & RARL (single-state) & \multicolumn{4}{c}{RARL (multi-state)} \\
\cmidrule(lr){2-2}\cmidrule(lr){3-6}
 & \textbf{Mean \boldmath$\pm$ SD}\tnote{a} & Seed 1 & Seed 2 & Seed 3 & \textbf{Mean \boldmath$\pm$ SD}\tnote{a} \\
\midrule
\multicolumn{6}{@{}l}{\textit{MC validation\tnote{b}}} \\
Feasibility rate & $(6.49 \pm 2.49)\%$ & $100.00\%$ & $100.00\%$ & $100.00\%$ & $100.00\%$ \\
Mean maneuver cost & $10.3320 \pm 0.0420$ & $10.3362$ & $10.3091$ & $10.2831$ & $10.3095 \pm 0.0266$ \\
MC cost spread, $\mathrm{SD}_{\mathrm{MC}}$\tnote{c} & $0.0829 \pm 0.0115$ & $0.1235$ & $0.1114$ & $0.1337$ & $0.1229 \pm 0.0112$ \\
95th-percentile cost & $10.4868 \pm 0.0443$ & $10.5506$ & $10.4978$ & $10.5074$ & $10.5186 \pm 0.0281$ \\
Central 90\% width\tnote{d} & $0.2549 \pm 0.0290$ & $0.3747$ & $0.3496$ & $0.4168$ & $0.3803 \pm 0.0340$ \\
\midrule
\multicolumn{6}{@{}l}{\textit{Nominal departure, $\Delta t_{\mathrm{dep}}=0$}} \\
Total maneuver cost & $10.2318 \pm 0.0218$ & $10.3154$ & $10.2976$ & $10.2707$ & $10.2946 \pm 0.0225$ \\
Cost above RARL (single-state)\tnote{e} & $0.0000$ & $0.0836$ & $0.0658$ & $0.0389$ & $0.0628 \pm 0.0225$ \\
Gap to RARL (single-state)\tnote{e} & $0.00\%$ & $0.82\%$ & $0.64\%$ & $0.38\%$ & $(0.61 \pm 0.22)\%$ \\
\bottomrule
\end{tabular*}
\begin{tablenotes}[flushleft]
\small
\setlength{\itemsep}{0pt}
\setlength{\parsep}{0pt}
\item[a] Aggregate entries give the mean $\pm$ sample SD of each row quantity across three independently trained policies.
\item[b] MC statistics use the same $10^4$ states per policy and include outcomes with impulse-cap violations.
\item[c] $\mathrm{SD}_{\mathrm{MC}}$ denotes the sample standard deviation of maneuver cost across MC states.
\item[d] Central 90\% width is the difference between the 95th and 5th cost percentiles.
\item[e] Differences use the single-state mean nominal cost as a fixed reference, shown as zero in its aggregate column.
\end{tablenotes}
\end{threeparttable}
\end{table}

All three multi-state RARL policies complete every MC transfer without truncation or impulse-cap violations, giving a feasibility rate of $100\%$ for each policy. This agreement across independent training runs supports the ability of RARL to learn reusable policies that retain feasibility beyond the departure states used for training and policy selection.
\par
\newpage
\noindent On the other hand, the single-state policies have feasibility rates of $3.75\%$, $8.62\%$, and $7.10\%$ for Seeds 1--3, respectively, with an aggregate value of $(6.49\pm2.49)\%$, as reported in Table~\ref{tab:em_multi_state_results}. Their infeasible outcomes arise from impulse-cap violations in the terminal departure or arrival impulse, with no Lambert failures or unbound-state truncations. Maneuver costs are therefore available for every MC outcome, and the reported cost statistics include the transfers with impulse-cap violations.

The multi-state policies have a mean MC maneuver cost of $10.3095\pm0.0266\,\mathrm{km/s}$. Each statistic is calculated over the MC states for one policy and then summarized across the three independent training runs. In particular, $\mathrm{SD}_{\mathrm{MC}}$ measures variation across departures, while the SD in the aggregate columns measures variation in the corresponding statistic across trained policies. The multi-state cost distributions have a larger $\mathrm{SD}_{\mathrm{MC}}$ and central 90\% width than their single-state counterparts. However, most single-state outcomes exceed the available terminal impulse authority, and their narrower cost distributions therefore cannot, on their own, be interpreted as more consistent feasible performance.

\begin{figure}[htbp]
  \centering
  \captionsetup{width=4.510in}
  \includegraphics[width=3.800in,trim=0 0 0 16bp,clip]{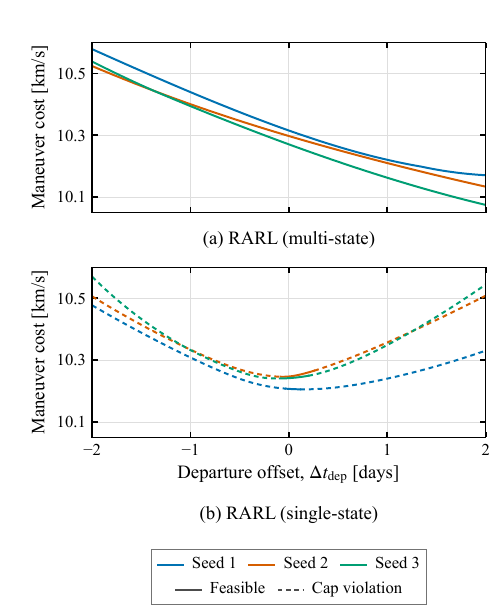}
  \caption{MC cost comparison of single- and multi-state RARL seeds.}
  \label{fig:res_cost_offset}
\end{figure}

Figure~\ref{fig:res_cost_offset} shows the dependence of maneuver cost on departure offset. For the multi-state policies, maneuver cost generally decreases toward later departure offsets within the tested interval, while feasibility is retained across the sampled departures. The lower single-state costs near nominal departure are consistent with specialization to the nominal training state. These policies also retain feasibility for some nearby departures, indicating limited local generalization, although most of the sampled departure family remains infeasible. The separation between solid and dashed portions of the curves in Fig.~\ref{fig:res_cost_offset} makes this distinction explicit.

At nominal departure, all three multi-state policies remain feasible, with a total maneuver cost of $10.2946\pm0.0225\,\mathrm{km/s}$. Relative to the single-state mean nominal cost of $10.2318\,\mathrm{km/s}$, the increase in the mean is $0.0628\,\mathrm{km/s}$, or $0.61\%$. For this case study, the broader departure-state feasibility is therefore accompanied by a small increase in mean nominal maneuver cost, illustrating an empirical robustness--cost tradeoff between nominal specialization and performance across dispersed initial states.

\begin{figure}[htbp]
  \centering
  \includegraphics[width=\linewidth,trim=0 0 0 14bp,clip]{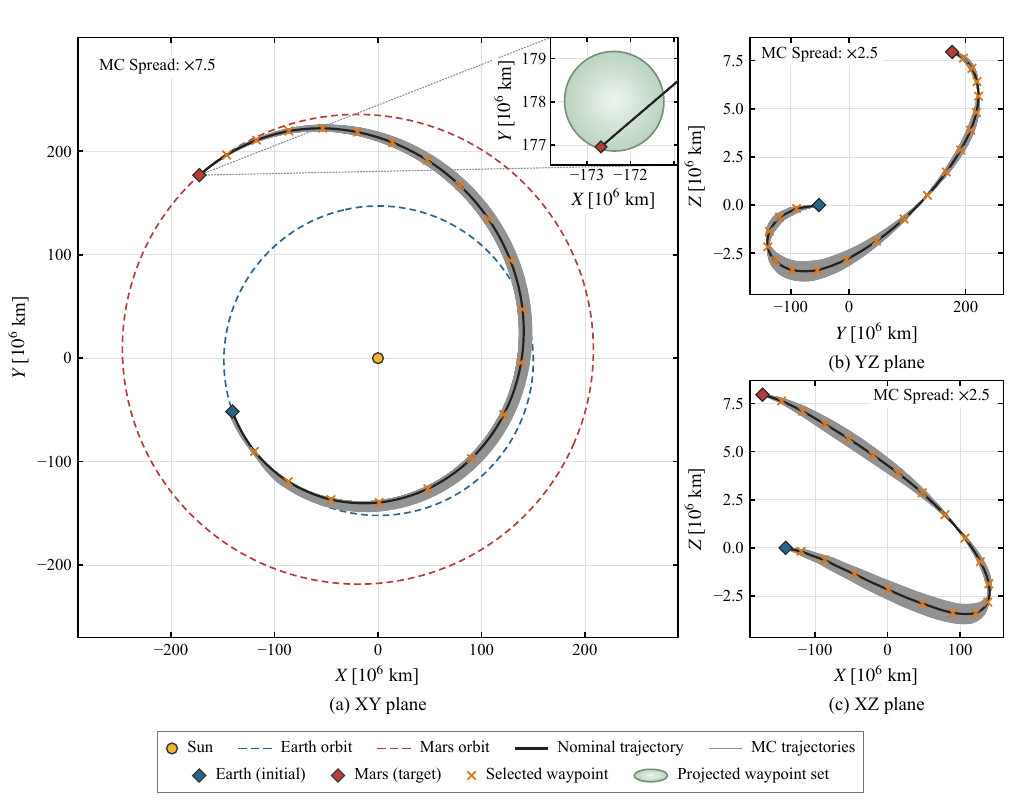}
  \caption{Earth--Mars trajectory ensemble for the representative multi-state RARL policy, showing 201 MC departures and nominal waypoints. MC deviations from the nominal trajectory are visually magnified $\times 7.5$ in the XY plane (a) and $\times 2.5$ in the YZ and XZ projections (b)--(c). The terminal inset shows the nominal waypoint set and MC trajectories at true scale.}
  \label{fig:res_traj_multi}
\end{figure}

\pagebreak

Figure~\ref{fig:res_traj_multi} illustrates the transfer family for Seed 2, the representative multi-state policy with the median MC mean maneuver cost. The figure shows trajectories from 201 sampled MC departures together with the nominal trajectory, while the waypoint markers and projected waypoint sets correspond to the nominal transfer. The trajectories reach the common terminal target from dispersed initial states, complementing the feasibility assessment over the full MC set in Table~\ref{tab:em_multi_state_results}.

\begin{figure}[htbp]
  \centering
  \includegraphics[width=\linewidth]{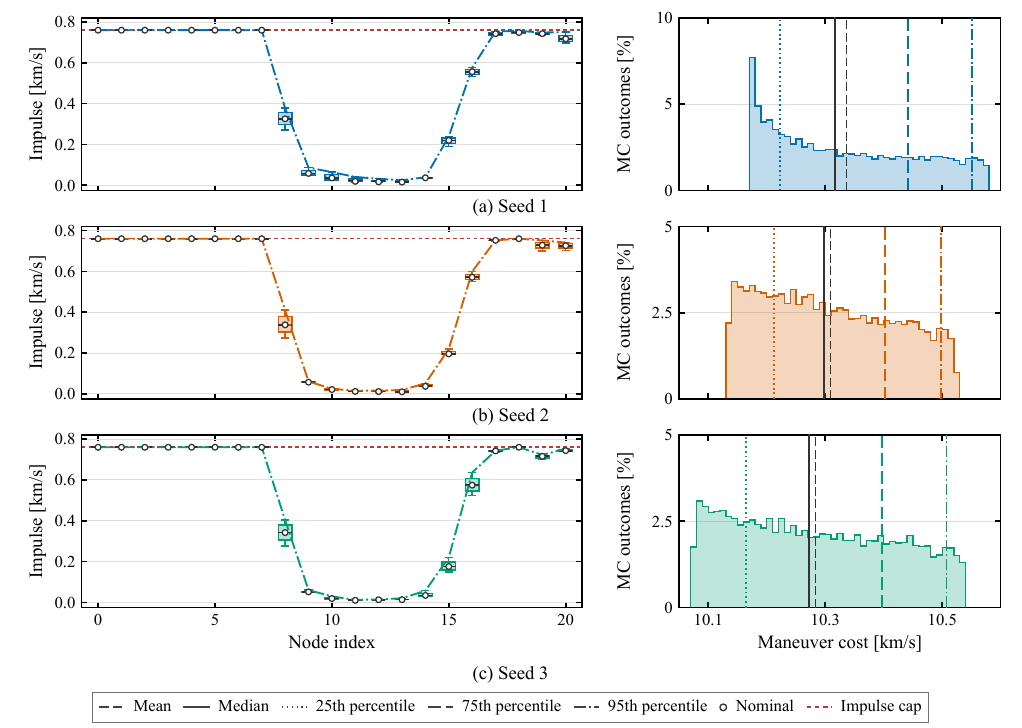}
  \caption{MC impulse and cost distributions for RARL (multi-state) policy seeds, with 5th--95th percentile whiskers and $0.01\,\mathrm{km/s}$ bins.}
  \label{fig:res_ctrl_multi}
\end{figure}

Figure~\ref{fig:res_ctrl_multi} pairs the node-wise impulse distributions with the maneuver-cost distribution for each multi-state policy. The allocation retains the approximate bang--off--bang structure observed at nominal departure. Impulses at nodes 0--7 remain close to the physical cap, followed by a reduction in maneuver effort. All three policies have median impulse magnitudes below the near-coast threshold of $0.05\,\mathrm{km/s}$ defined in Table~\ref{tab:em_single_state_results} at nodes 10--14. Maneuver effort then increases through nodes 15--18, and the terminal departure and arrival impulses, at nodes 19 and 20, remain below the cap for every MC realization. The boxes and percentile whiskers describe variation across departure states, while the nominal markers locate the corresponding nominal impulse within each distribution.

Successful terminal reconstruction alone does not establish feasibility, since the required departure or arrival impulse may exceed its cap. Across the sampled departures, the multi-state policies reach terminal conditions that admit the reconstructed two-impulse maneuver within both limits.

\FloatBarrier

Across three independent training runs, multi-state RARL extends feasible performance across the sampled departure family while retaining the approximate bang--off--bang maneuver structure and incurring a small increase in nominal cost.

%% file: sections/conclusions.tex
\section{Conclusion}
\label{sec:conclusions}

This paper presented Reachability Analysis-Informed Reinforcement Learning (RARL) for deterministic multi-impulse interplanetary transfers. The formulation places intermediate waypoint selection at the center of the learned decision process, using a local first-order reachability map to define the waypoint set and Lambert solutions to reconstruct the corresponding nonlinear transfer segments. This gives the policy a geometrically interpretable decision space while retaining classical astrodynamics for maneuver reconstruction. The accompanying diagnostics examine the action geometry, local approximation error, and terminal maneuver-demand assessment, making the role and limitations of this interface explicit.

In the Earth--Mars case study, this formulation recovers benchmark-quality nominal performance across three independent training runs, with a mean maneuver cost $1.72\%$ above the SCP reference. Training over dispersed initial states extends this capability beyond the nominal departure condition. Each of the three independently trained multi-state policies completes all $10^4$ Monte Carlo departures without truncation or impulse-cap violations, using the same selected policy across departures without further training. This broader sampled feasibility is accompanied by a $0.61\%$ increase in mean nominal maneuver cost relative to single-state training, while retaining the approximate bang--off--bang allocation. The resulting empirical robustness--cost tradeoff shows how the same reachability-informed formulation can support both nominal specialization and policy reuse across a family of departure conditions.

These findings concern a two-body Earth--Mars transfer with launch-time-induced initial-state dispersion, a fixed target state, and a fixed transfer duration. Within this setting, nonlinear feasibility is established empirically, and the reachability map remains a local approximation. Extending the formulation to broader transfer families, adaptive node placement, and refined reachability approximations would help assess its applicability to longer or more nonlinear arcs. The reusable state-feedback mapping also offers potential use as an onboard high-level guidance layer coupled to deterministic reconstruction, or as a structured initializer for subsequent trajectory optimization. Both applications build on the central contribution demonstrated here, a policy interface that couples learned selection of the transfer geometry with deterministic reconstruction of the maneuvers required to realize it.

%% file: sections/app-secondary-proxy.tex
\appendix
\section{Secondary Full-State Terminal-Difficulty Proxy}
\label{app:secondary_proxy}

The two-impulse construction in Section~\ref{subsec:mdp_reward} separates terminal position and velocity correction in a manner consistent with the terminal rendezvous. A supplementary alternative is obtained by estimating a single impulsive correction from the complete linearized terminal mismatch. Define
\begin{equation}
  \boldsymbol{G}_{k\rightarrow f}
  =
  \begin{bmatrix}
    \boldsymbol{\Phi}_{rv,k\rightarrow f} \\
    \boldsymbol{\Phi}_{vv,k\rightarrow f}
  \end{bmatrix}
  \in\mathbb{R}^{6\times3}.
  \label{eq:full_state_proxy_map}
\end{equation}
For this supplementary diagnostic, the matrix and state mismatch are evaluated using kilometers and seconds as the base units. This fixed unit convention determines the relative weighting of the position and velocity residuals in the least-squares fit. The minimum-norm least-squares estimate and its normalized magnitude are
\begin{subequations}
\label{eq:full_state_proxy}
\begin{align}
  \Delta\boldsymbol{v}_{6,k}^{\mathrm{lin}}
  & =
  \boldsymbol{G}_{k\rightarrow f}^{\dagger}
  \Delta\boldsymbol{x}_{k\rightarrow f},
  \label{eq:full_state_proxy_impulse} \\
  \eta_{6,k}
  & =
  \frac{
    \left\lVert\Delta\boldsymbol{v}_{6,k}^{\mathrm{lin}}\right\rVert_2
  }{\Delta V_{\max}}.
  \label{eq:full_state_proxy_metric}
\end{align}
\end{subequations}
Unlike the two-impulse assessment, $\eta_{6,k}$ represents the combined full-state correction through a single linearized impulse rather than through separate departure and arrival impulses. It therefore provides a complementary measure of terminal correction demand. The present reward adopts the two-impulse assessment to remain aligned with the terminal rendezvous structure, while $\eta_{6,k}$ remains available for diagnostic studies or for formulations that prefer a single-correction surrogate.